\documentclass[11pt]{amsart}
\usepackage[english]{babel}
\usepackage[toc,page]{appendix}
\usepackage[margin=1in]{geometry}
\usepackage{hyperref}
\usepackage{cite}
\usepackage{enumitem}
\usepackage{amsmath,amssymb,amsfonts}
\usepackage{graphicx}
\usepackage{textcomp}
\usepackage{xcolor}
\usepackage{mathtools, amsfonts, amssymb}
\usepackage{url}

\usepackage[normalem]{ulem}
\usepackage{caption}
\usepackage{algorithm}
\usepackage{algpseudocode}
\usepackage{mleftright}
\usepackage{tikz-cd}
\usepackage[normalem]{ulem}
\usepackage{pgf}
\usepackage{multirow}
\usepackage{array}

\usepackage{tikz}
\usetikzlibrary{positioning,fit,arrows.meta,backgrounds}
\usepackage{subfigure}
\usetikzlibrary{arrows,automata}
\usepackage{booktabs}
\usepackage{amsthm}
\usepackage{color}
\usepackage{pgfplots}
\usepackage{pgfplotstable}
\usepgfplotslibrary{fillbetween}
\pgfplotsset{compat=newest}
\usepackage{thmtools}
\usepackage{thm-restate}
\usepackage{scalerel,stackengine}
\usepackage{comment}

\stackMath
\newcommand\reallywidehat[1]{%
	\savestack{\tmpbox}{\stretchto{%
			\scaleto{%
				\scalerel*[\widthof{\ensuremath{#1}}]{\kern-.6pt\bigwedge\kern-.6pt}%
				{\rule[-\textheight/2]{1ex}{\textheight}}%
			}{\textheight}%
		}{0.5ex}}%
	\stackon[1pt]{#1}{\tmpbox}%
}
\usepackage{soul}
\usetikzlibrary{calc}
\makeatletter
\newif\if@anonymize
\@anonymizetrue    

\if@anonymize
\newcommand{\highlight@DoHighlight}{
	\fill [outer sep = -15pt, inner sep = 0pt, color=black]
	($(begin highlight)+(0,8pt)$) rectangle ($(end highlight)+(0,-3pt)$) ;
}

\newcommand{\highlight@BeginHighlight}{
	\coordinate (begin highlight) at (0,0) ;
}

\newcommand{\highlight@EndHighlight}{
	\coordinate (end highlight) at (0,0) ;
}

\newdimen\highlight@previous
\newdimen\highlight@current
\newlength{\item@width}

\DeclareRobustCommand*\anonymize{%
	\SOUL@setup
	\def\SOUL@preamble{%
		\begin{tikzpicture}[overlay, remember picture]
			\highlight@BeginHighlight
			\highlight@EndHighlight
		\end{tikzpicture}%
	}%
	\def\SOUL@postamble{%
		\begin{tikzpicture}[overlay, remember picture]
			\highlight@EndHighlight
			\highlight@DoHighlight
		\end{tikzpicture}%
	}%
	\def\SOUL@everyhyphen{%
		\discretionary{%
			\SOUL@setkern\SOUL@hyphkern
			\SOUL@sethyphenchar
			\tikz[overlay, remember picture] \highlight@EndHighlight ;%
		}{%
		}{%
			\SOUL@setkern\SOUL@charkern
		}%
	}%
	\def\SOUL@everyexhyphen##1{%
		\SOUL@setkern\SOUL@hyphkern
		\settowidth{\item@width}{##1}%
		\makebox[\item@width]{}%
		\discretionary{%
			\tikz[overlay, remember picture] \highlight@EndHighlight ;%
		}{%
		}{%
			\SOUL@setkern\SOUL@charkern
		}%
	}%
	\def\SOUL@everysyllable{%
		\begin{tikzpicture}[overlay, remember picture]
			\path let \p0 = (begin highlight), \p1 = (0,0) in \pgfextra
			\global\highlight@previous=\y0
			\global\highlight@current =\y1
			\endpgfextra (0,0) ;
			\ifdim\highlight@current < \highlight@previous
			\highlight@DoHighlight
			\highlight@BeginHighlight
			\fi
		\end{tikzpicture}%
		\settowidth{\item@width}{\the\SOUL@syllable}%
		\makebox[\item@width]{}%
		\tikz[overlay, remember picture] \highlight@EndHighlight ;%
	}%
	\SOUL@
}
\else
\newcommand{\anonymize}[1]{#1}
\fi
\makeatother
\newtheorem{theorem}{Theorem}[section]
\newtheorem{proposition}[theorem]{Proposition}
\newtheorem{proposition/definition}[theorem]{Proposition/Definition}
\newtheorem{lemma}[theorem]{Lemma}
\newtheorem{corollary}[theorem]{Corollary}

\newtheorem{problem}[theorem]{Problem}

\newtheorem{definition}[theorem]{Definition}

\newtheorem{example}[theorem]{Example}

\newtheorem{remark}[theorem]{Remark}

\declaretheoremstyle[
spaceabove=0pt,
spacebelow=0pt,
bodyfont=\normalfont,
postheadspace=0em,
qed=\qedsymbol,
headpunct={},
headformat={}
]{withouthead}

\newcommand{\h}{{\scriptscriptstyle\mathsf{H}}}

\DeclareMathOperator{\tr}{tr}

\let\diag\undefined%
\DeclareMathOperator{\diag}{diag}

\let\supp\undefined%
\DeclareMathOperator*{\supp}{supp}
\DeclareFontFamily{U} {MnSymbolC}{}
\DeclareFontShape{U}{MnSymbolC}{m}{n}{
	<-6> MnSymbolC5
	<6-7> MnSymbolC6
	<7-8> MnSymbolC7
	<8-9> MnSymbolC8
	<9-10> MnSymbolC9
	<10-12> MnSymbolC10
	<12-> MnSymbolC12}{}
\DeclareFontShape{U}{MnSymbolC}{b}{n}{
	<-6> MnSymbolC-Bold5
	<6-7> MnSymbolC-Bold6
	<7-8> MnSymbolC-Bold7
	<8-9> MnSymbolC-Bold8
	<9-10> MnSymbolC-Bold9
	<10-12> MnSymbolC-Bold10
	<12-> MnSymbolC-Bold12}{}
\DeclareSymbolFont{MnSyC} {U} {MnSymbolC}{m}{n}
\DeclareMathSymbol{\plus}{\mathrel}{MnSyC}{20}
\newcommand{\p}{\plus}

\newcommand{\pp}{{\plus\plus}}

\newcommand{\CG}{ \mathbb{C}[G]}
\newcommand{\card}[1]{{\left| #1 \right|}}

\newcommand{\Idg}{1_G}

\newcommand{\irr}[1]{\operatorname{Irr}\left( #1\right)}
\newcommand{\red}[1]{#1}
\newcommand{\blue}[1]{#1}
\newcommand{\cyan}[1]{#1}
\newcommand{\rred}[1]{{#1}}
 \newcommand{\rrred}[1]{{#1}}

\newcommand{\Rmnum}[1]{\expandafter\@slowromancap\Romannumeral #1@}

\tikzset{
	module/.style={%
		draw, rounded corners,
		minimum width=#1,
		minimum height=7mm,
		font=\sffamily
	},
	module/.default=2cm,
	>=LaTeX
}

\begin{document}

	\title{Sparse sum of Hermitian squares in group algebras of finite groups}
	\date{\today}
	\author[]{Jianting Yang}
	\address{ State Key Laboratory of Mathematical Sciences, Academy of Mathematics and Systems Science, Chinese Academy of Sciences, Beijing 100190, China }
	\email{yangjianting@amss.ac.cn}
	
	\author[]{Ke Ye}
	\address{State Key Laboratory of Mathematical Sciences, Academy of Mathematics and Systems Science, Chinese Academy of Sciences, Beijing 100190, China }
	\email{keyk@amss.ac.cn}
	
	\author[]{Lihong Zhi}
	\address{State Key Laboratory of Mathematical Sciences, Academy of Mathematics and Systems Science, University of Chinese Academy of Sciences, Beijing 100190, China }
	\email{\text{lzhi@mmrc.iss.ac.cn} }

	\thanks{\rrred{Jianting Yang, Ke Ye and  Lihong  Zhi  are supported by the National Key R\&D Program of China 2023YFA1009401. Jianting Yang is supported by the Postdoctoral Fellowship Program of CPSF under Grant Number GZC20252039, and the Basic Science Center Program (No: 12288201) of the National Natural Science Foundation of China.}}
	\thanks{To appear in Foundations of Computational Mathematics.}
\begin{abstract}
Nonnegative elements in group algebras play a central role in harmonic analysis, operator algebras, and computational optimization. This paper investigates sparse \rred{sum-of-Hermitian-squares} (SOHS) representations of nonnegative elements in the group algebras of finite groups. We prove that the convex relaxation of the sparse SOHS problem admits a closed-form solution, namely the square root of the given element. Based on this result, we propose a thresholding hierarchy for approximating sparse SOHS representations. We analyze the error of this hierarchy with respect to two natural residuals and establish exponential decay rates. Remarkably, one error bound is independent of the group order, and the other is also group-size independent when the group is cyclic or dihedral. These results extend existing work on Fourier sums of squares for abelian groups to a broader class of finite groups and provide new algorithmic tools for sparse noncommutative positivity certificates.
\end{abstract}  
	\maketitle
    \section{Introduction}\label{sec:intro}

In this work, we study sparse sums of Hermitian squares (SOHS) in the group algebras of finite groups. Let $G$ be a finite group and $\CG$ its group algebra. By definition, every element $f \in \CG$ can be uniquely expressed as $f=\sum_{x \in G}c_x x$ for some $ c_x \in \mathbb{C}$. For clarity, we introduce the coefficient function  $C_f:G \to \mathbb{C}$, which assigns to each group element $x \in G$ its corresponding coefficient $c_x$, i.e., $C_f(x) = c_x$. In this notation, $f = \sum_{x \in G} C_f(x) x$. The group algebra $\CG$ is equipped with the natural operations of addition, convolution (as multiplication), and the $*$-involution, given by $f^* = \sum_{x \in G} \overline{C_f(x)} x^{-1}$.

An element $f \in \CG$ is \emph{nonnegative} if, for every $*$-representation 
$\pi: \CG \to   \mathbb{C}^{n \times n}$, the matrix $\pi(f)$ is positive semidefinite. Nonnegativity in $\CG$ is equivalent to the existence of an SOHS decomposition \red{(cf. \cite[Proposition~4]{Schmüdgen2009} and  \cite[Corollary~2]{MR3067294})}
\[
f = \sum_{i=1}^r g_i^* g_i, \qquad g_i \in \CG,
\]
and hence provides algebraic certificates for \blue{nonnegativity}.

For abelian groups, all irreducible representations are one-dimensional. 
In this case, $\CG$ is isomorphic to the algebra of functions  on \rred{the dual group $\widehat{G}$}, and nonnegative elements in $\mathbb{C}[G]$ correspond to pointwise nonnegative functions on $\widehat{G}$. SOHS representations then reduce to Fourier sums of squares (FSOS) for functions on finite abelian groups, which have been extensively studied in polynomial optimization, Boolean function analysis, and combinatorial optimization \cite{MR3555385,yang2022computing,slot2023sum}. 

For non-abelian groups, higher-dimensional representations \red{arise}, and the equivalence between nonnegative elements in group algebras and pointwise \red{nonnegative functions on groups no longer holds}. In this setting, the SOHS \red{decomposition} provides \red{an important} link between representation theory and noncommutative optimization theory. 
More generally, Schmüdgen studied SOHS in finitely generated unitary $*$-algebras beyond finite groups~\cite{Schmüdgen2009}, while Ozawa analyzed the case of free group algebras~\cite{MR3067294}. Klep and Schweighofer~\cite{KS08} further established an equivalence between the existence of certain noncommutative SOHS and Connes' embedding conjecture in operator algebras. Recently, noncommutative SOHS decompositions have found important applications in quantum computing\cite{pb2g-j9cw}. 	

	\subsection{Main problem and our contributions}\label{subsec:intro-1.1}

Our focus is on \emph{sparse nonnegative elements}, i.e., nonnegative elements $f \in \CG$ whose support
\[
\operatorname{supp}(f)\coloneq \{x \in G: C_f(x)\neq 0\}
\]
is much smaller than $G$. 
Such elements naturally arise in combinatorial optimization, 
where objective functions are defined on exponentially large domains but only involve polynomially many terms. 
A prototypical example is the MAX-CUT problem, where the objective function over the Boolean cube 
$G = \{\pm 1\}^n$ has support of size $O(n^2)$ despite $|G| = 2^n$.

    As the following example illustrates,  a nonnegative element $f$ may not admit a sparse exact SOHS,  but it can have a sparse SOHS if a small perturbation is allowed.  

	\begin{example}[Motivating example  \rrred{\cite{MR2015916}}]\label{ex:motivating}
    Let \red{$G = \mathbb{Z}_2^n$}.  In this case,  we may identify $\CG$ with the space of polynomial functions on $\{0,1\}^n$.  The nonnegativity (resp.  SOHS) of elements in $\mathbb{C}[G]$ corresponds to the nonnegativity (resp.  SOS) of polynomial functions on $\{0,1\}^n$.  By \cite[Theorem~1.1 \& 1.2]{blekherman2016sums},  the polynomial 
		\[f(x_1,\dots, x_n) =\left( \sum_{j=1}^n x_j - \left\lfloor \frac{n}{2} \right\rfloor \right)  \left( \sum_{j=1}^n x_j - \left\lfloor \frac{n}{2} \right\rfloor  -1  \right)\]
    is nonnegative on $\{0,1\}^n$,  but there is no finite family $\{g_i: \deg g_i < n/2 \}_{i \in I}$ of polynomials such that $f = \sum_{i\in I} g_i^2$ on $\{0,1\}^n$.
		
		In contrast,  we have 
		\begin{equation}\label{ex:motivating:eq1}
			f (x_1,\dots,  x_n) + \frac{1}{4} = \left( \sum_{j=1}^n x_j - \left\lfloor \frac{n}{2} \right\rfloor -\frac{1}{2}\right)^2.
		\end{equation}
		Since  $f$ is integer-valued on $\{0,1\}^n$,  the sparse SOHS with error in \eqref{ex:motivating:eq1} certifies the nonnegativity of $f$.
	\end{example}
This motivates the following problem:
       \begin{problem}[Sparse SOHS problem]\label{prob:FSOS}
Given a nonnegative $f \in \CG$, find an SOHS representation $f  = \sum_{i \in I} g_i^\ast g_i$ with  $\card{ \cup_{i \in I} \supp(g_i)}$ 
 as small as possible.
\end{problem}

The existence of sparse SOHS certificates is crucial both for understanding noncommutative \blue{nonnegativity} 
and for \blue{producing} efficient algorithms in large-scale optimization. 
While sparse SOS and FSOS results are known for abelian groups and low-degree polynomials, 
little is understood for non-abelian groups or for general nonnegative elements without degree restrictions. 

Problem~\ref{prob:FSOS} has received significant attention \red{for abelian groups}. Fawzi et al.~\cite{MR3555385,MR5048294} 
studied the sparsity of exact SOS for nonnegative functions on finite abelian groups, 
showing that degree-$d$ functions on $\mathbb{Z}_N$ admit logarithmically sparse representations, 
and quadratic functions on $B_m^n$ also admit sparse SOS forms. 
Slot and Laurent~\cite{slot2023sum} analyzed the tightness of SOS bounds for degree-$d$ polynomials 
over the Boolean cube $\{-1,1\}^n$ and the $q$-ary cube $\mathbb{Z}_q^n$. Our work extends these results by considering arbitrary finite groups (including non-abelian groups) and removing degree constraints on $f$.

\subsection*{Our Contributions}  
The main contributions of this paper are summarized as follows:  

\begin{enumerate}[label=(\roman*)]
  \item \textbf{Convex relaxation and exact solution.}  
  We first formulate Problem~\ref{prob:FSOS} as the optimization problem~\eqref{Prob_min_SDP_Gram_modified}. In Theorem~\ref{Thm:square root}, we prove that the square root of $f$ is a closed-form exact solution to its convex relaxation, providing a noncommutative generalization of \cite[Proposition~3.4 \& Theorem~3.5]{yang2022computing}.  

  \item \textbf{Thresholding hierarchy.}  
  \red{Based} on Theorem~\ref{Thm:square root}, we propose a square-root-based hierarchy for approximately solving Problem~\ref{prob:FSOS} (Section~\ref{sec:SOHS-hierarchy}). This hierarchy employs thresholding techniques inspired by compressed sensing, \red{producing sparse SOHS decompositions  with increasing accuracy  when raising the order of the hierarchy.}

  \item \textbf{Error analysis.}  
  For $f \in \mathbb{C}[G]$ with $\Idg \succeq f \succeq \alpha \Idg$ for $\alpha \in [0,1)$, we establish two types of error estimates for the hierarchy in Theorems~\ref{thm:SDP_rate} and \ref{corr:error-bound-general}. When $\alpha>0$, both bounds decay exponentially with the relaxation order.  Remarkably, the bound in Theorem~\ref{thm:SDP_rate} is independent of the group order $|G|$, while the bound in Theorem~\ref{corr:error-bound-general} is also independent of $|G|$ when $G$ is cyclic or dihedral.
\end{enumerate}  

When $G$ is abelian, the condition $1_G\succeq f\succeq \alpha 1_G$ means that $f$, viewed as a function on the dual group, is pointwise bounded \rred{above by} $1$ and \rred{below} by $\alpha$. This assumption does not impose any sparsity condition on $f$, and analogous assumptions have appeared in previous work on Fourier sums of squares over finite abelian groups \cite{FF21,slot2023sum}. Our approach differs from the existing literature \rrred{in the following aspects}. First, previous works focus on finite abelian groups, and in some cases on specific families such as $\mathbb Z_d^n$ \cite{MR5048294,slot2023sum}, whereas our results apply to arbitrary finite groups. Second, existing hierarchies are typically degree-based, so their error bounds are expressed in terms of $\deg(f)$. By contrast, the square-root-based hierarchy used here leads to bounds depending on spectral information and, in Theorem~\ref{corr:error-bound-general}, on $|\operatorname{supp}(f)|$ rather than directly on the degree. Since $|\operatorname{supp}(f)|$ can be bounded in terms of $\deg(f)$, this also yields degree-based estimates, while the converse implication is not generally available. Third, besides the scalar residual $\delta_d$, which has been extensively studied in the literature, we also consider the coefficient $\ell_2$-residual $\varrho_d$. Theorem~\ref{thm:SDP_rate} shows that $\varrho_d$ admits a bound depending only on $d$ and the spectral lower bound $\alpha$, and is independent of both $|\operatorname{supp}(f)|$ and $\deg(f)$.

	\subsection{Related works}\label{sec:related:works}
	The study of the sum of squares (SOS) dates back to Hilbert's work  in 1888 \cite{hilbert1988}.  The core problem of SOS is concerned with the equivalence between two classes of real polynomials (resp. rational functions): nonnegative polynomials (resp. rational functions) and sums of squares of polynomials (resp. rational functions).  The well-known Motzkin polynomial disproved the equivalence for the polynomial case \cite{motzkin1967arithmetic}.  On the other hand,  the equivalence for the rational function case,  known as Hilbert's 17th problem, was proved by Artin in his seminal work \cite{artin1927zerlegung}. Based on various versions of \rrred{the Positivstellensatz}~\cite{Stengle1974,Schmudgen1991}, SOS has become a central approach in polynomial optimization in recent years \cite{Lasserre2001,Parrilo2003,scheiderer2009positivity,Laurent2009,Nie2014,klep2025sums}.

	As a generalization of SOS, the sum of Hermitian squares (SOHS) has been extensively studied \red{over recent decades}. It plays an essential role in various fields,  including complex geometry \cite{PS10}, noncommutative algebra \cite{Lewis88},  operator algebra \cite{KS08}, optimization \cite{KP10,GdL19}, and quantum information \cite{GdLM18,FF21,Bach25}. For general $\ast$-algebra $\mathcal{A}$, an SOHS of a nonnegative element  $f\in \mathcal{A}$ is defined as
	\[f = g_1^\ast g_1 + \cdots + g_r^\ast g_r,\]
	where $g_1,\dots,  g_r  \in \mathcal{A}$.  Arguably, the most studied cases are SOHS for polynomials \cite{yang2022computing,BR23} and noncommutative polynomials \cite{CKP12,GdLM18,KMP22}.  The former deals with the positivity of polynomial functions,  whereas the latter involves the positivity of operators.
   
    Group algebras, which form an important class of $\ast$-algebras, have been extensively studied in various fields of mathematics. \rred{We highlight four related directions} for \blue{SOHS} in group algebras. 

	First, when $G$ is finite abelian, the Fourier transform identifies $\CG$ with the algebra $F(\widehat G)$ of functions on the dual group $\widehat G$, equipped with pointwise multiplication. Under this identification, nonnegative elements in  $\CG$ correspond exactly to pointwise nonnegative functions on $\widehat G$. Moreover, after choosing a non-canonical isomorphism $G\cong\widehat G$, one may further regard these as functions on $G$. The basic differences between the abelian and non-abelian cases are summarized in Table~\ref{tab:compare}. Consequently, the SOHS \rred{representations} of nonnegative elements correspond to the Fourier sum of squares (FSOS) of pointwise nonnegative functions. Exact FSOS decomposition for pointwise nonnegative functions on finite abelian groups has been widely discussed in the literature.  On $\mathbb{Z}_2^n$, FSOS is equivalent to the SOS of polynomials on the binary hypercube $\{0,1\}^n$ \cite{MR3555385,blekherman2016sums}. Fawzi et al. established the degree bound for the exact FSOS of quadratic functions on $\mathbb{Z}_2^n$ \cite{MR3555385}, while {\rred{Sakaue}} et al. further extended this result to the degree bound for the exact FSOS for polynomials of arbitrary degree defined on $\mathbb{Z}_2^n$ \cite{MR3629437}. More generally, \rred{FSOS on $\mathbb{Z}_q^n$ has been studied in numerous works}. {\rred{Slot and Laurent}} analyzed the tightness of bounds provided by SOS \rred{hierarchies}  on $\mathbb{Z}_q^n$ \cite{slot2023sum}, and recently, Al-Sulami et al. proved the  existence of sparse exact FSOS for quadratic functions on $\mathbb{Z}_q^n$ \cite{MR5048294}. Fawzi et al. further showed that FSOS on finite abelian groups is closely related to chordal covers of Cayley graphs\cite{MR3555385}.  On the other hand, approximate FSOS has been investigated very recently. An efficient numerical algorithm for FSOS is proposed in \cite{yang2022computing}.

\begin{table}[ht]
		\centering
		\small
		\setlength{\tabcolsep}{5pt}
		\renewcommand{\arraystretch}{1.18}
		\caption{Positivity and SOS in Abelian and Non-abelian Group Algebras}
		\label{tab:compare}
		\begin{tabular}{@{}>{\raggedright\arraybackslash}p{0.23\textwidth}>{\raggedright\arraybackslash}p{0.37\textwidth}>{\raggedright\arraybackslash}p{0.28\textwidth}@{}}
			\toprule
			\textbf{group type} &  \textbf{abelian} & \textbf{non-abelian} \\ 
			\midrule
			algebra of functions vs. group algebra & 
			$F(G) \cong \mathbb{C}[G] $  (by Fourier transform and an identification $G \cong \widehat{G}$)&
			$F(G) \not \cong \mathbb{C}[G]$ \\ 
			\midrule
			characterization of nonnegativity & 
			\begin{enumerate}[label = (\alph*), leftmargin=*, nosep, topsep=0pt] 
				\item  $f$ (regarded as element in $F(G)$) is pointwise nonnegative.
				\item $f $ is a nonnegative element in $\CG$
			\end{enumerate} & 
			\begin{enumerate}[label = (\alph*), leftmargin=*, nosep, topsep=0pt] 
				\item  $f $ is a nonnegative element in $\CG$
			\end{enumerate}  \\ 
			\midrule
			terminology for SOS &  
			FSOS \red{(Fourier Sum of Squares)}
			&
			SOHS \red{(Sum of Hermitian Squares)}
			\\ 
			\bottomrule			
		\end{tabular}
\end{table}	
	
	Second, if $G$ is non-abelian, the $\ast$-algebra isomorphism between $F(G)$ and $\CG$ no longer holds. As a result, nonnegative elements in $\CG$ cannot be viewed as pointwise nonnegative functions on $G$. In this case, elements of $\CG$ are usually studied from the point of view of operator algebras. The \rrred{Positivstellensatz} for nonnegative elements of  $\CG$ is established in \cite{PS76}. The SOHS in the noncommutative group algebra is applied in the \rred{field} of the quantum information theory \cite{Watts2024relaxationsexact}. It is also shown that the Connes' embedding conjecture is equivalent to a problem involving SOHS in some free group algebras \cite{KS08}. Moreover, the technique of SOHS is one of the key tools in the study of algebraic structures with involution \cite{AU20}. Lastly, SOHS certificates play an important role in probability theory and interacting particle systems. For instance, the proof of Aldous' spectral gap conjecture~\cite{MR2629990} relies on the Octopus inequality~\cite[Theorem~2.3]{MR2629990}, which can be proved using the nonnegativity of certain elements in the group algebra of symmetric groups.

Third, nonnegative elements and SOHS certificates in group algebras provide natural tools for estimating spectral gaps. Let $G$ be a group, and suppose that $\Delta \in \CG$ satisfies $\Delta=\Delta^*$, and $\rho(\Delta)\succeq 0$ for each $\ast$-representation $\rho$. If for some $c>0$ the element $\Delta^2-c\Delta$ admits an SOHS certificate, then $c$ gives a lower bound for the spectral gap of $\Delta$. Spectral gaps play a central role in the analysis of random walks on groups and homogeneous spaces; see, for example, \cite{MR871832,MR4574854}. In group theory, Kazhdan's property (T) can be characterized by the existence of a spectral gap for an associated group Laplacian~\cite{MR3427595}. Semidefinite programming methods have also been developed to certify such gaps and compute related invariants, such as Kazhdan constants~\cite{MR3359223}.

The last direction focuses on the \blue{relation} between the \blue{SOHS} of $f \in \CG$ and the positive definiteness (cf.  Definition~\ref{def:pdfunction} and Proposition~\ref{prop:SOHSbasement}) of $C_f$. Following the basic work in \cite{MR23243}, positive definite functions have received significant attention in the fields of algebra and functional analysis \cite {MR1807896}.

	\section{Preliminaries}\label{sec:prelim}
	
	\subsection*{Notations}
	For clarity, we fix some notations.
	\begin{itemize}
    
		\item For a finite set $S$,  $\card{S}$ denotes the cardinality of $S$. 
		\item Given a complex matrix $A \in \mathbb{C}^{n\times n}$, $A^{\h}$ denotes the conjugate transpose of $A$, and $\|A\|_{2}$ denotes the largest singular value of $A$.
        \item For a finite group $G$, we denote by $\mathbb{C}^{G \times G}$ and $\mathbb{C}^{G}$ the spaces of matrices and \rred{vectors} whose elements are indexed by elements in $G \times G$ and $G$, respectively. Here, we choose and fix an order on $G$. Given $A \in \mathbb{C}^{G \times G}$, we use $A(x,y)$ to denote the element of $A$ indexed by $(x,y) \in G \times G$. Similarly, for $v \in \mathbb{C}^{G}$, we denote by $v(x)$ the element of $v$ indexed by $x \in G$.
		\item For a  complex vector $v\in\mathbb{C}^{n} $, we denote $\|v\|_2 \coloneqq \sqrt{v^\h v}$.
		\item We use $\Idg$ to denote the identity element in a group $G$.
        \item For a Hermitian element $f=f^* \in \CG$ and $\alpha \in \mathbb{R}$, the notation $f \succeq \alpha \Idg$ means that, for any $\ast$-representation $\pi$, the smallest eigenvalue of $\pi(f)$ is at least $\alpha$. Similarly, the notation $\Idg \succeq f$ means that for any $\ast$-representation $\pi$, the largest eigenvalue of $\pi(f)$ is at most $1$.
        \item We denote by $\mathsf{H}^n_{\p}$ and $\mathsf{H}^n_{\pp}$ the spaces of all $n \times n$ Hermitian positive semidefinite and positive definite matrices, respectively.
	\end{itemize}	

	\subsection{Group algebra and its representation}
	\red{ We first provide} a brief overview of group algebras. Let $G$ be a finite group.  The \emph{group algebra} of $G$ is defined as $\CG \coloneqq \cyan{ \{ \sum_{x\in G} a_x x: a_x\in \mathbb{C} \} }$, equipped with multiplication.
	\[
	m: \CG \times \CG \to \CG,\quad m\left( \sum_{x\in G} a_x x,  \sum_{x\in G} b_x x \right) \coloneqq \sum_{z\in G} \left( \sum_{xy = z} a_x b_y\right) z.
	\]
	For simplicity,  we denote $fg \coloneqq m(f,g)$ for $f,g\in \CG$.  \red{Formally, for every element $f \in \CG$, there exists a unique complex-valued function $C_f: G \to \mathbb{C}$,  such that $f$ can be explicitly formulated as} $f = \sum_{x \in G} C_f(x) x$. In the sequel,  we simply write $f = \sum_{x\in G} C_f(x) x$ for each $f \in \CG$.  We recall that $\CG$ is a $\ast$-algebra with the involution defined by $f^* \coloneqq \sum_{x \in G}\overline{ C_f(x)}x^{-1}$. 

	\begin{definition}[$\ast$-representation]\label{def:rep}
		A $\ast$-representation of $\CG$ is an algebra homomorphism $ \pi: \CG \to \mathbb{C}^{n \times n}$ such that \cyan{$\pi(x)$ is a unitary matrix for all $x \in G$} and $\pi(f^{\ast}) = \pi(f)^{\h}$, for each $f\in \mathbb{C}[G]$.
	\end{definition}
	Nonnegative elements in $\mathbb{C}[G]$ are defined in terms of $\ast$-representations.

	\begin{definition}[Nonnegative element]\label{def:nonnegative}
		An element $f \in \CG$ is nonnegative if for any $\ast$-representation $\pi: \mathbb{C}[G] \to \mathbb{C}^{n \times n}$,  $\pi(f)$ is positive semidefinite.  An SOHS representation of $f \in \CG$ is a decomposition $f=\sum_{i \in I} g_i^*g_i$ for some finite family $\{g_i\}_{i \in I} \subseteq  \CG$.
	\end{definition}

\begin{remark}
Finite abelian and non-abelian groups differ in SOHS decompositions.  For any $h \in \CG$, we define the  self-adjoint elements $\operatorname{Re}(h)=\frac{1}{2}\left(h+h^*\right)$ and $\operatorname{Im}(h)=\frac{1}{2i}\left(h-h^*\right)$. Then we have  
 $h=\operatorname{Re}(h)+i\operatorname{Im}(h)$. For any finite abelian group $G$  and $f \in \CG$ with 
$f=\sum_{j\in J}|h_j|^2$, we have the following decomposition
\[f=\sum_{j\in J}|h_j|^2 =\sum_{j\in J}\left| \operatorname{Re}(h_j) +i\operatorname{Im}(h_j)\right|^2=\sum_{j\in J} \operatorname{Re}(h_j)^2+ \sum_{j\in J} \operatorname{Im}(h_j)^2, \]
which converts a \emph{complex} SOHS decomposition into a \emph{real} SOHS decomposition. However, for general non-abelian groups, this property fails. Let $G$ be a finite  non-abelian group. For $f \in \CG$ with $f= hh^*$,  one can verify that 
 \[h h^* = \mathrm{Re}(h)^2 + \mathrm{Im}(h)^2 + \tfrac12(hh^* - h^*h),\]
where the commutator term  $\frac{1}{2}\left(hh^* - h^*h\right)$ is non-zero in non-abelian cases.  
\end{remark}

	Let $\irr{G}$ be the set of equivalence classes of irreducible representations of $G$. For each $\rho \in \irr{G}$,  we denote by ${\rho}:\CG \to  \mathbb{C}^{n_{\rho} \times n_{\rho}}$ the representation of $\CG$. According to \cite[Proposition 10]{serre1977linear},  the map \begin{equation}\label{eq:Phi}
		\Phi \coloneqq \oplus_{\rho \in \irr{G}} {\rho}:\CG \to \oplus_{\rho \in \irr{G}} \mathbb{C}^{n_{\rho} \times n_{\rho}},\quad \Phi(f) \coloneqq  \oplus_{\rho \in \irr{G}} {\rho}(f)
	\end{equation}	
	is a $\ast$-isomorphism. In particular,  for any $x \in G$,  $\Phi(x)^\h =\Phi(x)^{-1}=\Phi(x^{-1})$.  
	
	We define a Hermitian inner product on $\CG$:
	\[
	\langle \cdot,  \cdot \rangle:  \CG \times \CG \to \mathbb{C},\quad  \left\langle f,g \right\rangle \coloneqq C_{f^{\ast} g}(\Idg)=\sum_{x \in G}\overline{ C_f(x) } C_g(x),
	\] 
	where $\Idg$ is the identity element of $G$.  By \cite[Theorem 3]{serre1977linear} we have  
	\begin{equation}\label{eq:inner product}
		\left\langle f,g \right\rangle=\frac{1}{\card{G}}\sum_{\rho \in \irr{G}} n_{\rho}  \tr\left({\rho}(f)^\h {\rho}(g)\right). 
	\end{equation}
	
	\begin{example} \label{example:g-in-D6}
		Let $D_6$ be the dihedral group of order $6$, generated by $\sigma$ and $\tau$ with relations $\sigma^3=\tau^2=(\sigma \tau)^2=1_{D_6}$. Then $D_6$ is non-abelian and the group algebra $\mathbb{C}[D_6]$ is a $6$-dimensional vector space whose elements are linear combinations of $1_{D_6},\sigma,\sigma^2,\tau,\sigma\tau,\sigma^2 \tau$. \rred{The group} $D_6$ has three irreducible representations, \rred{namely}
		\begin{enumerate}[label = (\alph*)]
			\item The trivial representation $\rho_{0}:\mathbb{C}[D_6] \to \mathbb{C}$ determined by $\rho_{0}(\sigma)= \rho_0(\tau) = 1$.
			\item The sign representation   $\rho_{1}:\mathbb{C}[D_6] \to \mathbb{C}$ determined by $ \rho_{1}(\sigma)=1$, $\rho_{1}(\tau)=-1$.
			\item The standard  representation   $\rho_{2}:\mathbb{C}[D_6] \to \mathbb{C}^{2 \times 2}$ determined by 
			\[ \rho_{2}(\sigma)=\left(\begin{array}{cc}
				-\frac{1}{2}  & -\frac{\sqrt{3}}{2}  \\
				\frac{\sqrt{3}}{2}  & -\frac{1}{2} 
			\end{array}\right),\quad \rho_{2}(\tau)=\left(\begin{array}{cc}
				1 & 0 \\
				0 & -1
			\end{array}\right).\]
		\end{enumerate}
		Thus, we have a $\ast$-isomorphism
		\begin{equation}\label{ex:eq1}
			\Phi: \mathbb{C}[D_6] \to \mathbb{C}\oplus\mathbb{C} \oplus \mathbb{C}^{2\times 2},\quad \cyan{ \Phi(f)\coloneqq  \rho_0(f)\oplus \rho_1(f)\oplus \rho_2(f).}		\end{equation}
		As an \rred{illustrative} example, we take $f = \frac{1}{5} 1_{D_6}+ \frac{1}{10}\sigma+\frac{1}{20} \tau+ \frac{1}{20}\sigma\tau + \frac{1}{20} \sigma^2\tau$, then 
		\[
		\Phi(f)= \frac{1}{20} \cdot \left( \begin{array}{cccc}
			9 & 0 & 0 & 0 \\
			0 & 	3 & 0 & 0\\
			0& 0 & 3&   -\sqrt{3}\\
			0& 0 & \sqrt{3} &  3\\
		\end{array}\right).
		\]
	\end{example}
	
	We end this section with the definition of positive definite functions, which are also of great importance in representation theory and functional analysis \cite{jorgensen2016extensions}.

	\begin{definition}[Positive definite function]\cite[Definition~1.2]{jorgensen2016extensions}\label{def:pdfunction}
		For any group $G$, a complex-valued map $p:G \to \mathbb{C}$ is a positive definite function if for any finite subset $X \coloneqq \{x_1,\dots,x_r\} \subseteq G$, the matrix $K_X \in \mathbb{C}^{r \times r}$ is Hermitian positive semidefinite, where $K_X(i,j)\coloneq p(x_i^{-1}x_j)$, $1 \le i, j \le r$.
	\end{definition}	

Given a finite group $G$ and a function $p:G \to \mathbb{C}$, we define an element $f \coloneqq \sum_{x\in G}p(x)x\in \CG$. Let $R_f: \CG \to \CG$ be the linear map defined by the right multiplication by $f$. 
Indeed, for any $h \in \CG$,
\[  R_f(h)=hf= \sum _{x,y \in G }C_h(x) C_f (x^{-1}y)\cdot y \in \CG.  \]
Then the matrix $K_G=\left(C_f(x^{-1}y)\right)_{x,y\in G}$ in Definition~\ref{def:pdfunction}, with $X = G$, is precisely the right multiplication matrix with respect to the basis $G$ of $\CG$. 
	\subsection{Nonnegative elements in $\CG$}
	
	Let $G$ be a finite group. We recall from Definition~\ref{def:nonnegative} that an element $f$ is said to be \textit{nonnegative} if $\pi(f) \succeq 0$ for any $\ast$-representation $\pi$. We have the following equivalent characterizations of the nonnegativity in $\CG$.

	\begin{proposition}\label{prop:SOHSbasement}
		Let $G$ be a finite group and $f=f^* \in \CG$. Then the following statements are equivalent:
		\begin{enumerate}[label = (\alph*)]
			\item \label{prop:SOHSbasement-item1} $f$ is nonnegative.
			\item \label{prop:SOHSbasement-item2}$\Phi(f)$ is Hermitian positive semidefinite.
			\item \label{prop:SOHSbasement-item3}$f$ has an SOHS representation in $\CG$, that is, there exists a finite family $\{g_i\}_{i \in I} \subseteq  \CG $ such that $f=\sum_{i \in I}g_i^*g_i$. 	Furthermore, when this condition holds, $f$ also has a rank one SOHS, i.e., $f=g^*g$ for some $g \in \CG$.
			\item  \label{PD-function-item} $C_f:G \to \mathbb{C}$ is a positive definite function.
		\end{enumerate}
		
	\end{proposition}
	\begin{proof}
		$\ref{prop:SOHSbasement-item1}\Rightarrow	\ref{prop:SOHSbasement-item2}$ is clear. To prove $\ref{prop:SOHSbasement-item2}\Rightarrow	\ref{prop:SOHSbasement-item3}$, we observe that $\Phi(f) \succeq 0$, then there exists a matrix  $A \in \oplus_{\rho \in \irr{G}} \mathbb{C}^{n_{\rho} \times n_{\rho}}$ such that $\Phi(f)=A^\h A$. Since $\Phi$ is a $\ast$-isomorphism, we have $f=g^*g$ where $g \coloneqq \Phi^{-1}(A)$. The implication $\ref{prop:SOHSbasement-item3}\Rightarrow	\ref{prop:SOHSbasement-item1}$ can be proved by noticing that $\pi(f)=\sum_{i \in I}\pi(g_i)^\h\pi(g_i)$ for any $\ast$-representation $\pi$.

To complete the proof, it suffices to prove the equivalence of \ref{prop:SOHSbasement-item1} and \ref{PD-function-item}. Let
\[
K_G:=\left(C_f(x^{-1}y)\right)_{x,y\in G}.
\]
By Definition~2.5, $C_f$ is positive definite if and only if
$K_G \succeq 0$.    For any $u \in \mathbb{C}^G$,  we define $h\coloneq \sum_{x \in G} \overline{u(x)}\cdot x \in \CG$. As $K_G$ is the matrix corresponding to  right multiplication,  
\[  \sum _{x,y \in G }\overline{u(x)} K_G (x,y)\cdot y= \sum _{x,y\in G }C_f(x^{-1}y)\overline{u(x)}\cdot y =hf \in \CG.  \]

Therefore 
\begin{align*}
    u^{\h}K_Gu=\sum_{x \in G}  C_{h f}(x) u(x) &= C_{ h f h^*}(\Idg)\\
    &=\frac{1}{|G|}
\sum_{\rho\in\operatorname{Irr}(G)}
n_\rho\tr\left(\rho(h f h^* )\right)\\
&=\frac{1}{|G|}
\sum_{\rho\in\operatorname{Irr}(G)}
n_\rho\tr\left(
\rho( h)\rho(f)\rho( h)^{\mathrm H}
\right).
\end{align*}

To conclude,  $K_G \succeq 0$ if and only if  $\tr\left(
\rho( h)\rho(f)\rho( h)^{\mathrm H}
\right) \geq 0$ holds for all $h \in \CG$. Since the  direct sum of all irreducible representations, i.e., $\Phi $  is   a $\ast$-isomorphism, 
the condition that $\tr\left(
\rho( h)\rho(f)\rho( h)^{\mathrm H}
\right) \geq 0$ for all $h \in \CG$  
is equivalent to $\rho(f) \succeq 0$ for all  $\rho \in \irr{G}$.
	\end{proof}
    
	\begin{remark}
		 When $G$ is abelian, we recall that the nonnegativity of elements in $\mathbb{C}[G]$ is equivalent to the nonnegativity of functions on $\widehat{G}$, or on $G$ by identification $G \cong \widehat{G}$. \cyan{Keeping this in mind}, Proposition~\ref{prop:SOHSbasement} is a generalization of \cite[Proposition~3]{MR3555385} for non-abelian \cyan{groups}.
	\end{remark}	
	\begin{example}
		For illustrative purposes, we consider $f = \frac{1}{5} 1_{D_6}+ \frac{1}{10}\sigma+\frac{1}{20} \tau+ \frac{1}{20}\sigma\tau + \frac{1}{20} \sigma^2\tau \in \mathbb{C}[D_6]$ in Example~\ref{example:g-in-D6}. Denote $g\coloneq f+f^\ast$. Then it is easy to verify by \eqref{ex:eq1} that $1_{D_6} \succeq \cyan{g}  \succeq \frac{3}{10}1_{D_6}$. By definition, $C_g: \cyan{D_6} \to \mathbb{C}$ is given by
		\[
		C_g(x) = 
		\begin{cases}
			\frac{2}{5},\quad &\text{if~}x = 1_{D_6} \\
			\frac{1}{10},\quad &\text{otherwise}
		\end{cases}
		\]
		The matrix $K_G \in \mathbb{C}^{\cyan{D_6\times D_6} }$ defined in the proof of Proposition~\ref{prop:SOHSbasement} is 	
        \[
        K_G=\frac{1}{10}\cdot
\left(\begin{array}{cccccc}
4 & 1 & 1 & 1 & 1 & 1 \\
1 & 4 & 1 & 1 & 1 & 1 \\
1 & 1 & 4 & 1 & 1 & 1 \\
1 & 1 & 1 & 4 & 1 & 1 \\
1 & 1 & 1 & 1 & 4 & 1 \\
1 & 1 & 1 & 1 & 1 & 4 
\end{array}\right)
        \]
		where columns and rows of $K_G$ are indexed by elements $1_{D_6},\sigma,\sigma^2,\tau,\sigma\tau,\sigma^2 \tau$ of $D_6$. Accordingly, the SOHS of $g$ given by $K_G$ is simply
		\begin{eqnarray*}
			g =\frac{1}{60}\left(\sum_{x \in D_6}x\right)^*\left(\sum_{x \in D_6}x\right)+\frac{3}{10}\cdot\rred{ 1_{D_6}}.
		\end{eqnarray*}  
		
\end{example}
\section{Sum of Hermitian squares over group algebras of finite groups}\label{sec:square-root}
We first introduce the definition of the Gram matrices. For a fixed $f \in \CG$, any matrix $Q \in \mathbb{C}^{G \times G}$ satisfying the conditions
\begin{align}\label{Prob_feas_SDP_Gram}
Q^\h = Q,\;
	Q \succeq 0	,\; \sum_{x^{-1}y= z}Q(x,y)= C_f(z), \; \forall z \in G,
\end{align}
is called a Gram matrix of $f$, where $Q \in \mathbb{C}^{G \times G}$ means $Q$ is a complex matrix $\card{G}\times \card{G}$ whose rows and columns are indexed by elements of $G$. From the computational perspective, determining whether $f \in \CG$ is nonnegative is equivalent to determining whether $f$ admits a Gram matrix.  In fact, according to Proposition~\ref{prop:SOHSbasement}, the nonnegativity of $f$ is equivalent to the existence of an SOHS representation of $f$. It is straightforward to verify that any SOHS representation $f=\sum_{i\in I} g^*_i g_i$ provides a Gram matrix $Q$ satisfying \eqref{Prob_feas_SDP_Gram} such that $Q(x,y) = \sum_{i\in I}\overline{ C_{g_i}(x)} C_{g_i}(y)$.  Conversely, suppose that $Q$ satisfies \eqref{Prob_feas_SDP_Gram}, and $Q = \sum_{i\in I} u_i^\ast u_i$ is a decomposition of $Q$ where $u_i \in \mathbb{C}^G$ is a row vector indexed by $G$ for each $i \in I$.  Then $g_i \coloneqq \sum_{x\in G} u_i(x) x$ satisfies $f=\sum_{i\in I} g^*_i g_i$.

On the other hand, for nonnegative $f \in \CG$, we may reformulate the sparse SOHS problem (cf. Problem~\ref{prob:FSOS}) as
\begin{align}\label{prob:FSOS for f}
	\operatorname{minimize} \quad &\lvert \cup_{i \in I} \supp(C_{g_i})  \rvert  \\ 
	\operatorname{subject~to} \quad  &f = \sum_{i \in I} g_i^\ast g_i.\nonumber
\end{align}
Suppose $f = \sum_{i \in I} g_i^\ast g_i$ and let $Q$ be the corresponding Gram matrix.  Then we have 
\begin{equation}\label{eq:suppGram}
	\cup_{i \in I} \supp(C_{g_i}) = \{x\in G:  Q(x,x) \ne 0\}.
\end{equation} 
Combining \eqref{prob:FSOS for f} and  \eqref{eq:suppGram},  we obtain the following characterization of the sparse SOHS problem. 
\begin{lemma}\label{lem:reformulation}
	For any nonnegative $f\in\CG$,  \eqref{prob:FSOS for f} is equivalent to 
	\begin{align}\label{Prob_min_SDP_Gram}
		\operatorname{minimize} \quad &\|\diag(Q)\|_0 \\
		\operatorname{subject~to} \quad &Q \in \mathbb{C}^{G \times G},  \;
		Q^\h = Q,\;			Q \succeq 0 \nonumber \\
		\quad   &\sum_{ x^{-1}y=z}Q(x,y)= C_{f}(z),\;  z \in G,  \nonumber
	\end{align}
	where $\diag(X)$ denotes the vector consisting of diagonal elements of a matrix $X$ and $\lVert u \rVert_0 \coloneqq |\{1 \le j \le m: u_j \ne 0\}|$ is the $\ell_0$-norm of a vector $u\in \mathbb{C}^m$.  
\end{lemma}

\begin{remark}\label{rmk:trivial relaxation}
	We observe that the objective function in \eqref{Prob_min_SDP_Gram} involves the $\ell_0$-norm,  which is non-continuous and non-convex.  Solving such a problem is notoriously known to be challenging. Following the standard methods	\cite{MR2230846}, one may relax \eqref{Prob_min_SDP_Gram} by replacing $\ell_0$-norm with  $\ell_1$-norm.  However,  in our case,  the objective function of this relaxation is a constant function,  since $\lVert \diag(Q)  \rVert_1 = \sum_{x\in G}Q(x,x)=  C_{f}(\Idg)$.  As a consequence,  every feasible point of~\eqref{Prob_min_SDP_Gram} becomes an optimal solution under the $\ell_1$-norm relaxation,  yielding no meaningful distinction between them.
\end{remark}	

\begin{proposition}[Reformulation of sparse SOHS problem]
	For any nonnegative $f \in \CG \setminus \{0\}$,  the optimal value of 
	\begin{align}\label{Prob_min_SDP_Gram_modified}
		\operatorname{minimize} \quad &\card{ \{x \in G: Q(x,x)>0,~x\neq \Idg \}}+1 \\
		\nonumber	\operatorname{subject~to} \quad &Q \in \mathbb{C}^{G \times G},  \; Q^\h = Q,\; Q \succeq 0,~ Q(\Idg ,\Idg) >  0, \\  
		\nonumber	&\sum_{x^{-1}y=z}Q(x,y)= C_f(z), \; z \in G,
	\end{align}
	is the same as that of \eqref{prob:FSOS for f}.
\end{proposition}
\begin{proof}
	By Lemma~\ref{lem:reformulation}, it is sufficient to prove that \eqref{Prob_min_SDP_Gram} and \eqref{Prob_min_SDP_Gram_modified} have the same optimal value.  We denote the optimal value of \eqref{Prob_min_SDP_Gram} (resp.  \eqref{Prob_min_SDP_Gram_modified})  by $e$ (resp.  $e'$),  which is achieved by $Q$ (resp.  $Q'$).  Then we have 
	\[
	e' = \card{ \{x \in G: Q'(x,x)>0,~x\neq 1_G \}}+1 = \lVert \diag(Q') \rVert_{0} \ge e, 
	\]
	since $Q'$ is also a feasible point of \eqref{Prob_min_SDP_Gram}.
	
	If $Q(1_G,  1_G) \ne 0$,  then $Q$ is a feasible point of \eqref{Prob_min_SDP_Gram_modified},  which implies \rred{$e \ge e'$}.  Thus,  we may assume that $Q(1_G,  1_G) = 0$.  Since $f \ne 0$,  we have $Q \ne 0$.  Thus,  there is some $y_0 \in G$ such that $Q(y_0, y_0) > 0$ as $Q \succeq 0$.  We consider the matrix $\widetilde{Q} \in \mathbb{C}^{G \times G}$ defined by
	\begin{equation}\label{eq:Qtilde}
		\widetilde{Q}(x,y) = Q(y_0x,  y_0 y),\quad x,y\in G.
	\end{equation}
	By definition,  $\widetilde{Q}$ is the matrix obtained by simultaneously permuting rows and columns of $Q$.  Therefore,  we have $\widetilde{Q}^\h = \widetilde{Q} \succeq 0$ and $\widetilde{Q}(1_G,  1_G) = Q(y_0, y_0) > 0$. Moreover,  for any $z \in G$, 
	\[\sum_{x^{-1}y=z}\widetilde{Q}(x,y)= \sum_{x^{-1}y= z}Q(y_0  x,y_0y)= \sum_{x^{-1}y=z}Q(x,  y) = C_f(z).\]
	This implies that $\widetilde{Q}$ is a feasible point of \eqref{Prob_min_SDP_Gram_modified} such that 
	\[
	e' \le \card{ \{x \in G: \widetilde{Q}(x,x)>0,~x\neq 1_G \} }+1 = \|\diag(\widetilde{Q})\|_0 =   \|\diag(Q)\|_0 = e.  \qedhere
	\]
\end{proof}
We note that the objective function of \eqref{Prob_min_SDP_Gram_modified} is still the $\ell_0$-norm of a vector.  Directly solving \eqref{Prob_min_SDP_Gram_modified} is notoriously challenging.  Consequently,  we consider its $\ell_1$-norm convex relaxation:  	
\begin{align}\label{Prob_min_SDP_Gram_relaxed0}
	\operatorname{minimize} \quad  &\|\diag(Q)\|_1-Q(\Idg,\Idg) \\
	\operatorname{subject~to}\quad  &Q \in \mathbb{C}^{G \times G},  \; Q^\h = Q,\; Q \succeq 0,~ Q(\Idg ,\Idg) > 0,  \nonumber	\\
	&\sum_{x^{-1}y= z}Q(x,y)= C_f(z), \; z \in G. 			\nonumber
\end{align}

Since $\|\diag(Q)\|_1=\sum_{x\in G}Q(x,x)=  C_f(\Idg)$ is a constant,  \eqref{Prob_min_SDP_Gram_relaxed0} is equivalent to 
\begin{align} 	\label{Prob_min_SDP_Gram_relaxed}
	\operatorname{maximize} \quad &Q( 1_G,1_G) \\	
	\nonumber	\operatorname{subject~to} \quad &Q \in \mathbb{C}^{G \times G},  \;
	Q^\h = Q,\; Q \succeq 0,  \\  
	\nonumber	&\sum_{x^{-1}y=z}Q(x,y)= C_f(z), \; z \in G
\end{align}
We notice that if the optimal value of \eqref{Prob_min_SDP_Gram_relaxed} is $0$,  then we must have $f = 0$. Indeed,  given a fixed feasible point $Q$ of \eqref{Prob_min_SDP_Gram_relaxed} and $y_0 \in G$,  $\widetilde{Q}$ defined by \eqref{eq:Qtilde} is also a feasible point of \eqref{Prob_min_SDP_Gram_relaxed},  and $ Q(y_0,  y_0) = \widetilde{Q}(\Idg,  \Idg) = 0$.  Letting $y_0$ \cyan{run} through $G$,  we conclude that $Q = 0$ and $f = 0$.  As a consequence,  the constraint $Q(1_G,  1_G) > 0$ in \eqref{Prob_min_SDP_Gram_relaxed0} does not appear in  \eqref{Prob_min_SDP_Gram_relaxed}. 
\begin{remark}
	By Remark~\ref{rmk:trivial relaxation}, the $\ell_1$-norm relaxation of \eqref{Prob_min_SDP_Gram} is trivial, since $\lVert \operatorname{diag}(Q) \rVert_1$ is a constant. However, after reformulating \eqref{Prob_min_SDP_Gram} as \eqref{Prob_min_SDP_Gram_modified} by removing the term $Q(1_G,1_G)$, we obtain a non-trivial $\ell_1$-norm relaxation \eqref{Prob_min_SDP_Gram_relaxed0}, or equivalently, \eqref{Prob_min_SDP_Gram_relaxed}. In this reformulated relaxation, the objective is to minimize $C_f(\Idg)-Q(\Idg,\Idg)$, or equivalently to maximize $Q(\Idg,\Idg)$, since $C_f(\Idg)$ is a constant once $f$ is fixed,  and $Q(\Idg,\Idg)$ is not fixed by the linear constraints.
\end{remark}

The remainder of this section focuses on proving Theorem~\ref{Thm:square root}, which provides a closed-form solution to \eqref{Prob_min_SDP_Gram_relaxed}. To accomplish this,  we first recall the following fact.\footnote{We thank Sizhuo Yan for her great help.} 
\begin{theorem}[Alberti’s theorem] \cite{alberti1983note} \label{thm:Alberti-theorem}
	Given $P,Q \in \mathsf{H}^n_{\p}$,  we have 
	\[
	\tr\left( \sqrt{\sqrt{P} Q\sqrt{P}}\right) = 
	\sqrt{\inf_{ Y\in \mathsf{H}^n_{\pp} } \tr(YP) \tr(Y^{-1}Q)}.
	\]
\end{theorem}
Consequently, we derive the lemma that follows.
\begin{lemma}\label{lem:Yans-Lemma}
	Let $n, r$ be positive integers and let $A_1,\dots,  A_r,  B \in \mathsf{H}^n_{\p}$.  If $B^2 = \sum_{j=1}^r A_j^2$,  then 
	\[\tr(B)^2 \geq \sum_{j=1}^r\tr(A_j)^2.\]
\end{lemma}
\begin{proof}
For each $1 \le j \le r$, we take $P= I_n$ and $Q=A_j^2$ in Theorem~\ref{thm:Alberti-theorem}. Then for any $S \in H^n_{\pp}$, we obtain 
\[ 
\tr(A_j)^2 = \inf_{\rred{Y} \in \mathsf{H}^n_{\pp}} \tr(Y) \tr(Y^{-1}A_j^2) \le \tr(S) \tr(S^{-1}A_j^2).
\]
This together with $B^2 = \sum_{j=1}^r A_j^2$ implies
\[ 
\sum_{j=1}^r \operatorname{tr}(A_j)^2 \leq \operatorname{tr}(S) \sum_{j=1}^r \operatorname{tr}(S^{-1}A_j^2) =  \operatorname{tr}(S) \operatorname{tr} (S^{-1}B^2).
\] 
Since $S \in H^n_{++}$ is arbitrary, we have 
\[
\sum_{j=1}^r \operatorname{tr}(A_j)^2 \le \inf_{S \in \mathsf{H}^n_{\pp}} \operatorname{tr}(S) \operatorname{tr} (S^{-1}B^2) = \tr(B)^2,
\]
where the equality follows from Theorem~\ref{thm:Alberti-theorem} by taking $P= I_n$ and $Q= B^2$. 
\end{proof}
We also recall the following trace inequality.
\begin{lemma}\cite[Corollary 2.1.24]{bikchentaevtrace}\label{lem:Abs-ineq}
	For any matrix $A \in \mathbb{C}^{n\times n} $,  we have  
	$\tr(|A|) \geq |\tr(A)|$,  where $|A|\coloneqq \sqrt{A^*A}$.
\end{lemma} 
With all the above preparations, we finally arrive at the promised theorem on the optimal solution of \eqref{Prob_min_SDP_Gram_relaxed}. 
\begin{theorem} \label{Thm:square root}Let $G$ be a finite group and  $f\in \CG$ be a nonnegative element.
	\begin{enumerate}[label = (\alph*)]
		\item There exists a unique nonnegative element $h \in \mathbb{C}[G]$ such that $h=h^\ast$ and $f=h^2$.  \label{Thm:square root:item1}
		\item The matrix $Q_h \in \mathbb{C}^{G \times G}$ is an optimal solution of \eqref{Prob_min_SDP_Gram_relaxed},  where $Q_h$ is defined by $Q_{h} (x,y)= C_{h} (x^{-1})C_h (y)$ for $x,y \in G$ and $h = \sum_{z\in G} C_h(z) z$.   \label{Thm:square root:item2}
	\end{enumerate}
\end{theorem}
\begin{proof}
	We first prove \ref{Thm:square root:item1}.  Let $\Phi$ be the $\ast$-isomorphism $\Phi:\CG \to \oplus_{\rho \in \irr{G}} \mathbb{C}^{n_\rho \times n_{\rho}} $ defined in~\eqref{eq:Phi}.  Since $f$ is nonnegative, by Proposition~\ref{prop:SOHSbasement}, the element $h \coloneqq \Phi^{-1}(\sqrt{\Phi(f)})$ satisfies the condition. \rred{The uniqueness follows from the uniqueness of the positive semidefinite square root of each block  $\Phi(f)_\rho$.}
		
	To prove \ref{Thm:square root:item2}, \rred{the case $f=0$ is trivial}.  We notice that $Q_h$ is a feasible point of \eqref{Prob_min_SDP_Gram_relaxed}:
	\[
	\sum_{x^{-1}y = z} Q_h(x,y) = \sum_{x^{-1} y = z} C_h(x^{-1}) C_h(y) = C_{h^2}(z) = C_f(z),\quad z \in G.
	\]
	Suppose $\widetilde{Q}$ is a feasible point of \eqref{Prob_min_SDP_Gram_relaxed} \cyan{with} rank $r \coloneqq \operatorname{rank}(\widetilde{Q}) \ge 1$.  Let $\widetilde{Q}=\sum_{i=1}^r  u_i u_i^{\h}$ be a spectral decomposition of $\widetilde{Q}$,  where $u_1,\dots,  u_r$ are nonzero column vectors in $\mathbb{C}^G$.  Given $x \in G$ and $v \in \mathbb{C}^G$,  we denote by $v(x)$ the element of $v$ labeled by $x$.  Therefore,  we have $\widetilde{Q}(1_G,1_G)=\sum_{i=1}^r \left| u_i(1_G)\right|^2$.  We consider 
	\[
	g_i \coloneqq \sum_{x\in G} \overline{u_i(x)} x \in \mathbb{C}[G],\quad 1 \le i \le r. 
	\]
	By definition,  we have $C_{g_i} (x) = \overline{u_i(x)}$.  Constraints of \eqref{Prob_min_SDP_Gram_relaxed} indicate that $f=\sum_{i=1}^r g_i^{\ast} g_i$.  Let $h \in \mathbb{C}[G]$ and $\Phi = \oplus_{\rho \in \irr{G}} {\rho}$ be as in the proof of \ref{Thm:square root:item1}.  The construction of $h$ implies $\rho(h) \in \mathsf{H}^n_{\p}$.  Using $h^2 = f = \sum_{i=1}^r g_i^{\ast} g_i$ and  \eqref{eq:inner product},  we obtain for each $ \rho \in \irr{G}$ that
	\begin{align*}
		\rho (h)^2 =\sum_{i=1}^r\rho(g_i)^{\h} \rho(g_i),\quad 
		C_h(1_G) =\sum_{\rho \in \irr{G}} \frac{n_\rho}{\card{G}}\tr(\rho(h)). 
	\end{align*} 
	For each $\rho \in \irr{G}$,  we consider the embedding 
	\[
	\iota_\rho: \mathbb{C}^{n_{\rho} \times n_{\rho}} \to \mathbb{C}^{n_{\rho}^2 \times n_{\rho}^2},\quad  \iota_{\rho}(A) \coloneqq A \otimes I_{n_\rho}.
	\]
	Here $A \otimes B$ denotes the Kronecker product of two matrices $A$ and $B$.  Therefore,  we have 
	\[
	\Psi: \mathbb{C}[G] \to \oplus_{\rho \in \irr{G}} \mathbb{C}^{n_{\rho}^2 \times n_{ \rho}^2 },\quad \Psi(f) = \rred{\oplus}_{\rho \in \irr{G}} \iota_\rho  ( \rho (f) ) =  \rred{\oplus}_{\rho \in \irr{G}}  \rho (f) \otimes I_{n_\rho}.
	\]
	It is straightforward to verify that $\Psi$ is a $\ast$-homomorphism and 
	\begin{equation}\label{Thm:square root:eq1}
		\tr(\Psi(f))= \sum_{\rho \in \irr{G}} n_\rho\tr(\rho(f)).  
	\end{equation}
	Thus,  we may derive
	\begin{align}\label{Thm:square root:eq2}
		\widetilde{Q}(1_G,1_G)=\sum_{i=1}^r \left| u_i(1_G)\right|^2 =\sum_{i=1}^r \left| C_{g_i}(1_G) \right|^2 &=\sum_{i=1}^r  \left| \sum_{\rho \in \irr{G}} \frac{n_\rho}{\card{G}}\tr(\rho(g_i)) \right|^2\\
		&=\frac{1}{\card{G}^2} \sum_{i=1}^r \left|   \tr(\Psi(g_i) )\right|^2 \nonumber \\
		& \leq \frac{1}{\card{G}^2} \sum_{i=1}^r   \tr( \left|\Psi(g_i)\right| )^2,  \nonumber
	\end{align} 
	where the last inequality follows from Lemma~\ref{lem:Abs-ineq}.  Moreover,  the equation $h^2=\sum_{i=1}^r g_i^{\ast} g_i$ leads to 
	\[ \Psi(h)^2=\sum_{i=1}^r \Psi(g_i)^\h \Psi(g_i)=\sum_{i=1}^r\left|\Psi(g_i)\right|^2. \] 
	According to Lemma~\ref{lem:Yans-Lemma},  we conclude that $\tr( \Psi(h)  )^2 \geq \sum_{i=1}^r   \tr( \left|\Psi(g_i)\right| )^2$,  which combined with \eqref{Thm:square root:eq1}, \eqref{Thm:square root:eq2} and Lemma~\ref{lem:Abs-ineq} implies $\widetilde{Q}(1_G,1_G) \leq  {Q}_{h}(1_G, 1_G)$. \qedhere
\end{proof}   
\begin{remark}
	If $G$ is a finite abelian group, Theorem~\ref{Thm:square root} reduces to \cite[Theorem~3.5]{yang2022computing}, whose proof relies on the Fourier analysis of functions on finite abelian groups. The existence of SOHS representations in Theorem~3.8(a) may be viewed as a finite-group-algebra analogue of the operator Fejér--Riesz theorem \cite{Dritschel2010,MR1815959,IJS25}. However, it is not a direct consequence of the operator version. Indeed, the operator Fejér--Riesz theorem concerns univariate Laurent polynomials with operator coefficients, whereas our results concern elements of finite group algebras. Moreover, the optimality statement in Theorem~\ref{Thm:square root}(b) is a distinctive feature of the finite-group-algebra setting and has no direct counterpart in the operator Fejér--Riesz theorem.
\end{remark} 

\section{Thresholding \rred{hierarchies} for   sparse SOHS \rred{problems}}\label{sec:SOHS-hierarchy}
For convenience,  in the rest of the paper,  we denote $\sqrt{f} \coloneqq h$,  where $h$ is the unique element constructed in Theorem \ref{Thm:square root}. Although we have already derived a closed-form solution to the convex relaxation problem of~\eqref{prob:FSOS for f}, it still faces the following \rred{challenges}:
\begin{enumerate}[label = (\roman*)]
	\item  \textbf{High  computational complexity}: To solve the problem \eqref{prob:FSOS for f}, we need to explicitly compute~$\sqrt{f}$. \red{
However, in the sparse-support setting (where the number of non-zero coefficients of $f$ is much smaller than $|G|$), computational complexity scales rapidly with $|G|$, making it prohibitively expensive when $|G|$ is exponentially large. For example, when $G$ is an abelian group, the computational complexity is quasi-linear \rred{in} $\card{G}$.}
	
	\item \textbf{Non-sparsity}:
	Although $f$ may be extremely sparse, the support of $\sqrt{f}$ can nonetheless be large.
\end{enumerate}

In the following example over $\mathbb{Z}_N$, even when $f$ has only a constant number of nonzero terms, its square root $\sqrt{f}$ can have a support size of at least $\rred{\Omega}(\sqrt{N})$.

\begin{example}
\rred{Let $N$ be sufficiently large and set  $G=\mathbb{Z}_N=\{0,1,2,...,N-1\}$}, the characters of $G$ are given by $\chi_k(x)=\exp\left(\frac{2\pi k i x}{N}\right)$ for $k=0,1,2,...,N-1$. Consider the following group element $f \in \mathbb{C}[\widehat{G}]$ with 
$$f(x):=\frac{1}{2}+\frac{1}{4}\chi_1(x)+\frac{1}{4}\chi_{-1}(x)=\frac{1}{2}+\frac{1}{2}\cos\left(\frac{2\pi x}{N}\right)=\left|\cos\left(\frac{\pi x}{N}\right)\right|^2, ~x \in \mathbb{Z}_N.$$
Note that the sparsity of $f$ is $3$. 
Since
 \begin{equation*}
	\left\vert \cos (x)\right\vert^2 =\cos^2 (x)
=\frac{1+\cos(2x)}{2},
\end{equation*}
we have
$$\sqrt{f}(x) =\left|\cos\left(\frac{\pi x}{N}\right)\right|.$$
By the Fourier series of the $|\cos(t)|$ on the interval $[0,2\pi]$
 \begin{equation*}
	\left\vert \cos(t)\right\vert =\frac{2}{\pi }+\frac{4}{\pi }
	\sum_{k=1}^{\infty }\frac{(-1)^{k}}{1-4k^{2}}\cos (2kt),~ t \in [0,2\pi],
\end{equation*}
for $x \in \mathbb{Z}_N$, we have 
 \begin{eqnarray*}
 	\left|\cos\left(\frac{\pi x}{N}\right)\right| &=&\frac{2}{\pi }+\frac{4}{\pi }
 	\sum_{k=1}^{\infty }\frac{(-1)^{k}}{1-4k^{2}}\cos \left(\frac{2k\pi x}{N}\right)\\
 	&=&\frac{2}{\pi }+\frac{4}{\pi }
 	\sum_{k=1}^{N-1 }\frac{(-1)^{k}}{1-4k^{2}}\cos\left(\frac{2k\pi x}{N}\right)
 	+\frac{4}{\pi }\sum_{k=N}^{\infty }\frac{(-1)^{k}}{1-4k^{2}}\cos \left(\frac{2k\pi x}{N}\right)\\
 	 &=&\frac{2}{\pi }+\frac{4}{\pi }
 	\sum_{k=1}^{N -1}\frac{(-1)^{k}}{1-4k^{2}} \left(\frac{\chi_k(x)}{2}+\frac{\chi_{-k}(x)}{2}\right)
 	+\frac{4}{\pi }\sum_{k=N}^{\infty }\frac{(-1)^{k}}{1-4k^{2}}\cos \left(\frac{2k\pi x}{N}\right).
 \end{eqnarray*}
 The last equation holds because $\frac{1}{2}\chi_{k}(x)+\frac{1}{2}\chi_{-k}(x)=\cos\left(\frac{2k\pi x}{N}\right)$. 
For $k \in \{1,2,3,...,\lfloor  \frac{N-1}{2} \rfloor\}$  the Fourier coefficient $\widehat{\sqrt{f}}(\chi_k)$ of $\sqrt{f}$ at  $\chi_{k}$ satisfies:
 \begin{eqnarray*}
\rred{ 	\widehat{\sqrt{f}}(\chi_k)= \frac{4}{\pi}\cdot \frac{(-1)^{k}}{2-8k^{2}}+\frac{4}{\pi}\cdot \frac{(-1)^{N-k}}{2-8(N-k)^2}+\frac{4}{\pi}\sum_{\substack{m \geq  N \\ m \equiv \pm k \pmod{N}}} \frac{(-1)^{m}}{2-8m^{2}}.}
 \end{eqnarray*}

 Next, we prove that if $N \geq 3$ and $\rred{k < \sqrt{\frac{3N-4}{8}}}$, then  $\widehat{\sqrt{f}}(\chi_k) \neq 0$. Let
 \begin{eqnarray*}
 A_k=\frac{4}{\pi}\cdot \frac{(-1)^{k}}{2-8k^{2}}, ~B_k=\rred{\frac{4}{\pi}\cdot \frac{(-1)^{N-k}}{2-8(N-k)^2}+\frac{4}{\pi}\sum_{\substack{m \geq  N \\ m \equiv \pm k \pmod{N}}} \frac{(-1)^{m}}{2-8m^{2}}}
 \end{eqnarray*}
then $\widehat{\sqrt{f}}(\chi_k)= A_k+B_k$. Obviously, if $|A_k| >|B_k|$, then $\widehat{\sqrt{f}}(\chi_k)= A_k+B_k \neq 0$. For $k \leq N/2$, define 
 \[   r_N\coloneq \frac{4}{\pi} \sum_{m=\rred{\lceil N/2\rceil}}^{\infty }\frac{1}{8m^{2}-2} \geq \rred{\frac{4}{\pi}\cdot \left|\frac{(-1)^{N-k}}{2-8(N-k)^2}\right|}+ \frac{4}{\pi} \sum_{\substack{m \geq  N \\ m \equiv \pm k \pmod{N}}} \left| \frac{(-1)^{m}}{2-8m^{2}} \right| \geq |B_k|,\]
 we have the following estimation about $r_N$:
\[ r_N=\frac{4}{ \pi}\sum_{m=\rred{\lceil N/2\rceil}}^{\infty }\frac{1}{8m^{2}-2} \leq \frac{4}{ \pi} \sum_{m=\rred{\lceil N/2\rceil}}^{\infty } \frac{1}{6m^2}  \leq \frac{4}{ \pi} \int_{\rred{\lceil N/2\rceil}-1}^{\infty}\frac{1}{6m^2} \leq \frac{4}{ \pi}  \cdot \frac{1}{3N-6}.\]

For any  $k < \sqrt{\frac{3N-4}{8}}$, we have  $|A_k| > r_N \geq |B_k|$, and thus $\widehat{\sqrt{f}}(\chi_k) \neq 0$. 

In conclusion, the sparsity of $f$ is $3$. In contrast, the sparsity of $\sqrt{f}$ is at least \rred{$\Omega(\sqrt{N})$}, indicating that computing $\sqrt{f}$  becomes computationally infeasible when $N$ is exponentially large.
\end{example}

Our approach to circumvent these limitations is inspired by compressed sensing. The thresholding method is a post-processing technique used in compressed sensing to sparsify dense solutions obtained by solving relaxation problems. The key idea is to retain only components of large magnitude of dense solutions, thereby enforcing  sparsity. For further details, interested readers are referred to \cite{MR2230846}.

In this section, we propose a hierarchy for \eqref{prob:FSOS for f} based on the thresholding method. The main idea is to define a sequence of subsets $\{S_d \subseteq G: d\in \mathbb{N}\}$ with the following properties: 
\begin{enumerate}[label = (\alph*)]
	\item $|S_d|$ is small for each $d\in \mathbb{N}$.
	\item \cyan{$S_d$} \rred{is} easy to compute.
	\item $S_d \subseteq S_{d+1}$ for each $d\in \mathbb{N}$ and \rred{$ \supp (\sqrt{f}) \subseteq \cup_{d \in \mathbb{N}} S_d $}.
\end{enumerate}
Subsequently, for each $d \in \mathbb{N}$, we solve the SDP feasibility problem \eqref{Prob_feas_sparse_SDP_Gram} below to obtain an SOHS supported in $S_d$, and this gives us a thresholding hierarchy for \eqref{prob:FSOS for f}.
\begin{align}\label{Prob_feas_sparse_SDP_Gram}
	\operatorname{minimize} \quad  &1\\
	\operatorname{subject~to}\quad  &Q\in \mathbb{C}^{G \times G},\;Q^\h = Q,\; Q \succeq 0,  \nonumber	\\
	&\sum_{x^{-1}y= z}Q(x,y)= C_f(z), \; z \in G, 			\nonumber\\
	& Q(x,x)=0,  \; x\notin S_d. \nonumber	
\end{align}

As shown in Theorem \ref{Thm:square root},  $\sqrt{f}$ is an optimal solution of \eqref{Prob_min_SDP_Gram_relaxed}. Our definition of $S_d$ is based on an approximation of $\sqrt{f}$. Without loss of generality, in the following discussion, we always assume that $1_G \succeq f$.

\begin{proposition}[Approximation of $\sqrt{f}$] \label{thm:error_of_coeffs}
	Suppose $G$ is a finite group and $0 \leq \alpha<1$ is a fixed real number.  Let $p$ be a univariate polynomial, and let $\varepsilon > 0$ be a real number such that $\left| \sqrt{\cyan{t}}-p(\cyan{t}) \right|\rred{\leq}  \varepsilon$ on $[\alpha,1]$.  Then for any nonnegative $f \in \CG$ such that  $\Idg \succeq f \succeq  \alpha \Idg $, we have
	\[
	\|\Phi(\sqrt{f})- \Phi (p(f))\|_2 \leq \varepsilon, 
	\]
	where $\Phi = \oplus_{\rho \in \irr{G}} {\rho}: \CG \to \oplus_{\rho \in \irr{G}} \mathbb{C}^{n_{\rho} \times n_{\rho}}$ is the isomorphism defined in \eqref{eq:Phi} and $n_{\Phi} \coloneqq \sum_{\rho \in \irr{G}} n_{\rho}$.  Moreover, $|C_{p(f)}(x)- C_{\sqrt{f}}(x)| \le \varepsilon$ for each $x \in G$.   
\end{proposition} 
\begin{proof}
	By assumption, $\Phi(f)$, $ \Phi(\sqrt{f})$, and $\Phi(p(f))$ are pairwise commuting Hermitian matrices.  If $\Phi(f) =U^{\h}\Lambda U$ is a spectral decomposition of $\Phi(f)$,  then $\Phi(\sqrt{f})=U^{\h} \sqrt{\Lambda} U$ (resp.  $\Phi(p(f))=U^{\h} p(\Lambda) U$) is a spectral decomposition of $\Phi(\sqrt{f})$ (resp.  $\Phi(p(f))$).  Thus,  we have  
	\[
	\| \Phi(\sqrt{f})-\Phi\left( p(f) \right)\|_2 
	=  \| U^\h\left(\sqrt{\Lambda}-p(\Lambda) \right)U \|_2
	= \|\sqrt{\Lambda}-p(\Lambda)  \|_2  \le \varepsilon.
	\]
	The last inequality holds since  $I_{n_{\Phi}} \succeq \Phi(f ) \succeq \alpha I_{n_{\Phi}}$.

	For \rred{the moreover statement}, we denote $r \coloneqq \sqrt{f}-p(f)$.  Clearly,  we have $ \| \Phi(r) \| \le \varepsilon$ and $\| \rho(r) \|_2 \le \varepsilon$ for each $\rho \in \irr{G}$.  Given $x\in G$,  by \eqref{eq:inner product} we have 
	\[
	\left| C_r(x) \right|=\left|\sum_{\rho \in \irr{G}} \frac{n_\rho}{\card{G}}\tr(\rho(r)^\h \rho(x)) \right|
	\leq \sum_{\rho \in \irr{G}} \frac{n_\rho}{\card{G}}n_\rho\|\rho(r)\|_2 \cdot \|\rho(x)\|_2
	= \sum_{\rho \in \irr{G}} \frac{n_\rho^2}{\card{G}} \| \rho(r) \|_2 
	\leq \varepsilon.    \]
	Here, the last inequality follows from the fact that $\rho(r)$ is a Hermitian matrix \blue{of size} $n_\rho \times n_\rho$ for each $\rho \in \irr{G}$.
\end{proof}

As a consequence of Proposition~\ref{thm:error_of_coeffs},  we obtain the following method to approximate $\sqrt{f}$ for $f \in \CG $ such that $\Idg  \succeq f \succeq \alpha \Idg$.  Let $\{p_d: \deg(p_d) \leq d\}_{d\in \mathbb{N}}$ be a sequence of univariate polynomials such that $|p_d(t) - \sqrt{t} | \le \varepsilon_d$ for any $t\in [\alpha,  1]$, where $\{\varepsilon_d\}_{d\in \mathbb{N}}$ is a decreasing sequence
of positive real numbers with $\lim_{d\to \infty} \varepsilon_d = 0$.  For each $d \in \mathbb{N}$,  we define 
\begin{equation}\label{equ:Base-select}
	S_d \coloneqq \bigcup_{k\leq d}\operatorname{supp}({p_k(f)} ) \subseteq G.
\end{equation}
Clearly,  the sequence $\{ S_d \}_{d\in \mathbb{N}}$ satisfies $S_d \subseteq S_{d+1}$ for each $d\in \mathbb{N}$. By Proposition~\ref{thm:error_of_coeffs},  $p_d(f)$ is a sequence in $ \CG$ approximating $\sqrt{f}$ in both the  operator norm of image under $\Phi$  and the $\ell_{\infty}$-norm of coefficients:
\[
\lim_{d\to \infty} \lVert \Phi(p_d(f))  -  \Phi(\sqrt{f}) \rVert_2 =  \lim_{d\to \infty} \max_{x\in G} \{ | C_{p_d(f)}(x) - C_{\sqrt{f}}(x) | \} = 0.
\] 
In the next section, we will quantitatively analyze the error of the thresholding hierarchy with respect to $S_d$ defined by \eqref{equ:Base-select}.

\begin{remark}
	In practice, to avoid computing $p_d$, we can simply take $\widetilde{S}_d \coloneqq \cup_{k\leq d}\operatorname{supp}(f^k) \subseteq G$, in which case it follows that~$	  \cup_{k\leq d} S_k \subseteq \cup_{k\leq d} \widetilde{S}_k $. It is also worth mentioning that, in the case where $\card{G}$ is large, the complexity of computing $S_d$ for small $d \in \mathbb{N}$  is much lower compared to computing $\sqrt{f}$.
\end{remark}

We conclude this section with a summary of the basic properties of $S_d$ defined by \eqref{equ:Base-select}.
\begin{proposition}\label{prop:S_d-thresholdingd}
	For $f \in \CG$ with $\Idg \succeq f \succeq \alpha \Idg$, we have
	\begin{enumerate}[label = (\alph*)]
		\item $|C_{\sqrt{f}}(x)|\rred{\leq} \varepsilon_d$ for any $x\notin S_d$. \label{prop:S_d-thresholdingd:item1}
		\item For any group isomorphism $\psi: G \to G'$, we have $\psi(p_d(f))=p_d(\psi(f))$, where  $\psi(f) \coloneqq \sum_{x\in G}C_{f}(x)\psi(x)$. \label{prop:S_d-thresholdingd:item2}
	\end{enumerate}
\end{proposition}
\begin{proof}
	\ref{prop:S_d-thresholdingd:item1} follows from $|C_{p_d(f)}(x)-C_{\sqrt{f}}(x)|\leq \varepsilon_d$ and \ref{prop:S_d-thresholdingd:item2} is clear from the definition. 
\end{proof}
\begin{remark}
	According to \ref{prop:S_d-thresholdingd:item1} of Proposition~\ref{prop:S_d-thresholdingd}, $S_d$ consists of components of a large magnitude in $\sqrt{f}$. This observation justifies the term ``thresholding hierarchy". Furthermore, \ref{prop:S_d-thresholdingd:item2} of Proposition~\ref{prop:S_d-thresholdingd} indicates that our thresholding hierarchy does not depend on how $G$ is presented.
\end{remark}	

\section{Error analysis}\label{sec:convergence analysis1}

In this section, we analyze the error of the thresholding hierarchy proposed in Section~\ref{sec:SOHS-hierarchy}. Since \eqref{Prob_feas_sparse_SDP_Gram} is an SDP feasibility problem, it is natural to measure the error of the hierarchy by the residuals of \eqref{Prob_feas_sparse_SDP_Gram}. Before we proceed, we define for each positive integer $d$ and each Hermitian positive semidefinite matrix $Q\in \mathbb{C}^{\cyan{S_d \times S_d}}$, the following nonnegative element:
\[
f_Q \coloneqq \sum_{x,y \in S_d} Q(x,y)x^{-1}y \in \mathbb{C}[G].
\]
Clearly, $f_Q$ has an SOHS representation supported on $S_d$. We consider two types of residuals:
\begin{itemize}
	\item The squared $\ell_2$-distance between $f$ and the feasible region of \eqref{Prob_feas_sparse_SDP_Gram}:
	\begin{eqnarray}\label{SDP:tau_d}
		\varrho_d \coloneq \min \left\{\sum_{x \in G} |C_{f- f_Q}(x)|^2 : Q\in  \mathbb{C}^{\cyan{S_d \times S_d}} ,\; Q=Q^\h,\; Q \succeq 0 \right\}.
	\end{eqnarray}
	Many SDP solvers and conic solvers \cite{sun2020sdpnal+,MR2895218} measure the infeasibility of primal problems in this way. Moreover, $\varrho_d$ can be regarded as a generalization of the $\ell_2$-error in compressed sensing. In Theorem~\ref{thm:SDP_rate}, we will prove that $\varrho_d$ decays exponentially and that the rate is independent of $G$.
	\item The distance between $f$ and the feasible region of \eqref{Prob_feas_sparse_SDP_Gram} along the constant term direction:
	\begin{eqnarray}\label{SOHS:M_d}
		\delta_d \coloneq \min \{ \lambda \ge 0: f+\lambda = f_Q,\; Q\in \mathbb{C}^{\cyan{S_d \times S_d}},\; Q=Q^\h,\ Q \succeq 0 \}.
	\end{eqnarray}
	In the convergence analysis of SOS hierarchies, the residual along the constant term direction is a widely used error measurement \cite{slot2023sum,BR23}. In Theorem~\ref{corr:error-bound-general},  we will show that $\delta_d$ decreases exponentially to zero. Furthermore, if $G$ is cyclic or dihedral, the convergence rate is independent of $G$ (cf. Corollary~\ref{corr:error-bound-ZN}).
\end{itemize}

\subsection{Polynomial approximation}
This subsection focuses on the approximation rate of some special functions by univariate polynomials, which plays an essential role in our analysis of $\varrho_d$ and $\delta_d$. To begin with, we recall two basic results in approximation theory. For $\alpha\in [0,1)$ and $d \in \mathbb{N}$, we define
\begin{align*}
	\gamma_{\alpha,d} &\coloneqq	\min_{q \in \mathbb{R}[t]_d  } \max_{t \in [-1,-\alpha] \cup [\alpha,1]}  \left| q(t)-   \left|t  \right|  \right|,  \\
	\tau_{\alpha,2d+1} &\coloneqq	\min_{p \in \mathbb{R}[t]_{2d+1}  } \max_{t \in [-1,-\alpha] \cup [\alpha,1]}  \left| p(t)- \operatorname{sign}(t)  \right|.
\end{align*}
Here $ \mathbb{R}[t]_{m} $ is the vector space of univariate polynomials with degree at most $m$.  \cyan{We have the  following results about  $\gamma_{0,d}$ and $\tau_{\alpha,2d+1}$:}

\begin{lemma}\cite{varga1987conjecture}\label{lem:Bernstein_constant}
The limit $\lim_{d\rightarrow \infty} d \gamma_{0,2d}$ exists and is non-zero.
\end{lemma}

\begin{lemma}\cite[Theorem~1]{eremenko2006uniform}\label{lem:Exp_rate_of_sign}
\cyan{	For any $0< \alpha <1$, we have}
	\[
	\lim_{d\rightarrow \infty} \sqrt{d} \left(\frac{1+\alpha}{1-\alpha}\right)^d \tau_{\alpha,2d+1} = \frac{1-\alpha}{\sqrt{\pi \alpha}}.
	\]
\end{lemma}
As a consequence, we obtain the following estimate of $\gamma_{\alpha,d}$.

\begin{lemma}[Approximation of $|t|$]\label{lem:app_abs_vals_alpha}
	For any $0 \le \alpha<1$, there exists a constant $c_{\alpha}>0$ such that for any $d\in \mathbb{N}$ \rred{with $d \geq 2$}, we have
	\[
	\gamma_{\alpha,d} \leq 
	\begin{cases}
		\frac{c_{\alpha}}{\sqrt{d-1}} \left(\sqrt{\frac{1-\alpha}{1+\alpha}}\right)^d \quad &\text{if $0 < \alpha < 1$},\\
		\frac{c_{\alpha}}{d} \quad &\text{if $\alpha = 0$}.\\
	\end{cases}
	\]
\end{lemma} 
\begin{proof}
	If $\alpha = 0$, the inequality is clearly justified by Lemma~\ref{lem:Bernstein_constant}. According to Lemma~\ref{lem:Exp_rate_of_sign},  there is a constant $c'_{\alpha} > 0$ such that for any $k\in \mathbb{N}$,  one can find a polynomial $u_{2k+1}(t)$ of degree at most $2k+1$ satisfying the inequality
	\[
	\left|\operatorname{sign}(t)-  u_{2k+1}(t) \right| \le c'_{\alpha} \frac{(1-\alpha)^{k+1}}{ (1 + \alpha)^k \sqrt{k \pi \alpha}},\quad t \in [-1,-\alpha]\cup[\alpha,1].   
	\]
	For each integer $d>0$, we let $k$ be the largest integer such that $2k + 2 \le d$.  Since $|t| =t \operatorname{sign}(t)$,  we have
	\[    
	\left| |t|-t\cdot u_{2k+1}(t) \right|   \leq |t|  \left|\operatorname{sign}(t)-  u_{2k+1}(t) \right| \le c'_{\alpha} \frac{(1-\alpha)^{k+1}}{ (1 + \alpha)^k \sqrt{k \pi \alpha}},\quad  t \in [-1,-\alpha]\cup[\alpha,1].
	\]
	The existence of $c_{\alpha}$ and the desired inequality follow immediately.   
\end{proof}
According to \cite{varga1987conjecture}, $c_0$ in Lemma~\ref{lem:app_abs_vals_alpha} actually lies in $(0.278,1/\pi)$. Next, we observe the approximation of the square root function.
\begin{lemma}\label{Lem:approx_square_root}
	Let $0\leq \alpha<1$ be a real number and let    $\{q_d\}_{d \in \mathbb{N}}$ (resp.  $\{r_d\}_{d \in \mathbb{N}}$)  be a sequence of univariate polynomials of degree at most $d$ (resp. positive real numbers) such that $||t|-q_d(t)|\leq r_d$ on $[-1,-\sqrt{\alpha}]\cup [\sqrt{\alpha},1]$.  There exists  a sequence of  univariate polynomials  $\{p_d\}_{d \in \mathbb{N}}$ of degree at most $d$ such that $\left|  \sqrt{t}-p_d(t) \right|\leq r_{2d}$ on $[\alpha,1]$.
\end{lemma}
\begin{proof}
	For each $d\in \mathbb{N}$,  without loss of generality, we can assume that $q_d$ is an even polynomial function, since $ (q_d(t)+ q_d(-t))/2$ is also a polynomial approximation of the absolute value function.
	Since $q_d(t)$ is an even polynomial function,  $p_d(t) \coloneqq q_{2d}(\sqrt{t})$ is a polynomial in $t$ and $p_d(t^2) = q_{2d}(t)$.  Given any $t \in [\sqrt{\alpha},1]$,  we have $\left| |t|- p_d(t^2)\right|<r_{2d}$,  from which we obtain $\left| \sqrt{t'}- p_d(t')\right|<r_{2d}$ for any $t' \coloneqq t^2 \in [\alpha,1]$. 
\end{proof}

Given $\alpha \in [0,1)$ and $d \in \mathbb{N}$, we define
\begin{eqnarray*} 
	E_{\alpha,d} \coloneq	\min_{p \in \mathbb{R}[t]_d  } \max_{t \in   [\alpha,1]}  \left| p(t)- \sqrt{t}  \right|.
\end{eqnarray*}
We also denote by $p_{\alpha,d} \in \mathbb{R}[t]_d$ the best approximation of $\sqrt{t}$ on $[\alpha,1]$. A combination of Lemmas \ref{lem:app_abs_vals_alpha} and \ref{Lem:approx_square_root} leads to the following estimate.

\begin{proposition}[Approximation of $\sqrt{t}$] \label{coro:appr-rate-sqrt}
	Given any $\alpha \in [0,1)$, there is a constant $c_\alpha > 0$ such that for each $d\in \mathbb{N}$, we have 
	\[
	E_{\alpha, d} \le \begin{cases}
		\frac{c_{\alpha}}{\sqrt{2d-1}} \left(\frac{1-\alpha^2}{1+\alpha^2}\right)^{d} \quad &\text{if $0 < \alpha < 1$}, \\
		\frac{c_{0}}{d} \quad &\text{if $\alpha = 0$}.\\
	\end{cases}
	\]
\end{proposition}

\begin{remark}
	An estimate for $E_{0,d}$ is also established in \cite{stahl2003best}. We also notice that the error bound in Proposition~\ref{coro:appr-rate-sqrt} can be improved for $\alpha> 1/5$.  Indeed,  for each positive integer $d$,  we let $q_d(t)$ be the degree $d$ Chebyshev interpolation polynomial for $\sqrt{t}$ on $[\alpha,1]$. According to \cite[Theorem 3.6]{10.5555/2161609}, we have
	\[			E_{\alpha,d} \le \left( \frac{1-\alpha}{2} \right)^{d +1} \frac{(2d -1)!!\alpha^{\frac{1}{2}- (d +1)}}{2^{2d +1} (d +1)!} 
	< \left( \frac{1-\alpha}{4\alpha} \right)^{d},
	\]
	where $n!! $ is the double factorial of $n$. 
	
	For independent interest, we mention that $E_{\alpha,d}$ can be numerically estimated by min-max optimization or semi-infinite programming \cite{GuoFengRecom2}. For comparison, we provide some numerical upper bounds of $E_{\alpha,d}$ in Table~\ref{tab:E_alpha_d_values}. \cyan{Moreover, $E_{\alpha,d}\leq 1/2$ for all $\alpha \in [0,1)$ and $d\in \mathbb{N}$, since $p(t)=1/2$ is a polynomial approximation of $\sqrt{t}$ on $[0,1]$  with uniform error at most $1/2$. }
	\begin{table}[ht]
		\centering
		\caption{Upper bounds of $E_{\alpha,d}$} \label{tab:E_alpha_d_values}
		\begin{tabular}{c | c c c c c}
			\toprule
			$E_{\alpha,d}$	& \multicolumn{5}{c}{$\alpha$} \\
			\midrule
			d & $0.0$ & $0.01$ & $0.1$& $0.2$ &$0.4$ \\
			\midrule
			$1$  &$0.126$&   $0.093$& $ 0.046$&$  0.028$&$  0.012 $\\
			$2$ & $0.069$&  $0.041$& $ 0.013$&  $0.007$& $0.003$\\
			$3$ & $0.047$& $ 0.023$& $0.005$& $0.003$& $0.002 $\\
			\bottomrule
		\end{tabular}
	\end{table}
\end{remark}

We notice that the key ingredient in the proof of Proposition~\ref{coro:appr-rate-sqrt} is the fact that one can approximate $\sqrt{t}$ on $[\alpha, 1]$ exponentially by a polynomial when $\alpha > 0$. Thus, it is tempting to expect an exponential approximation of $\sqrt{t}$ on $[0,1]$. Unfortunately, this is not possible, as shown by the following proposition.
\begin{proposition}[Impossibility theorem for uniform approximation]\label{thm:impossible_continuous}
	There exists no polynomial sequence $\{p_i: \operatorname{deg}(p_i) \leq i\}_{i\in \mathbb{N}}$ converging uniformly and exponentially to the square root function on the interval $[0,1]$. 
\end{proposition}

\begin{proof}
	Suppose there exists such a sequence of polynomials $\{q_d(t):\operatorname{deg}(q_d)\leq d\}_{d \in \mathbb{N}}$. Then the sequence
	$\{q_d(t^2)\}_{d \in \mathbb{N}}$ converges uniformly and exponentially to the absolute value function on the interval $[-1,1]$. This contradicts  Lemma~\ref{lem:Bernstein_constant}.
\end{proof}

\subsection{Convergence rate of $\varrho_d$}
In this subsection, we establish an estimate for the convergence rate of $\varrho_d$.
\begin{theorem}[Convergence rate of $\varrho_d$]\label{thm:SDP_rate}
	Given $\alpha \in [0,1)$, there is a constant $c_{\alpha}>0$ such that for any finite group $G$, $f\in \mathbb{C}[G]$ with $\Idg \succeq f \succeq \alpha \Idg$, and $d\in \mathbb{N}$, we have 
	\[
	\varrho_d \le 
	\begin{cases}
		\frac{c_\alpha}{{d}} \left(\frac{1-\alpha^2}{1+\alpha^2}\right)^{2d} \quad &\text{if $0 < \alpha <1$}, \\
		\frac{c_\alpha}{d^2} \quad &\text{if $\alpha = 0 $}. \\
	\end{cases}
	\]
\end{theorem}
\begin{proof}
	Let $h_d\coloneq p_{\alpha,d}(f)$. Then  $h_d=h_d^*$ and $h_d^2=|h_d|^2$. Thus we have,
	\[\|\Phi (f-h_d^2) \|_2 \leq \| \Phi(\sqrt{f}-h_d)  \|_2 \cdot  \|  \Phi(\sqrt{f}+h_d)  \|_2. \]
	As $p(t)=0$ is \cyan{a} uniform approximation of the square root \rred{function} with error at most $1$,  we know that $E_{{\alpha},d } <1$ for all $\alpha \in [0,1)$ and $d \in \mathbb{N}$. Since $ \| \Phi(\sqrt{f}-h_d)  \|_2 \leq  E_{{\alpha},d }$, and 
	$$\|  \Phi(\sqrt{f}+h_d)  \|_2\leq\| \Phi(\sqrt{f}) \|_2 + \|  \Phi(h_d) \|_2 \leq 2+E_{{\alpha},d},$$
	we can conclude that
	\[\| \Phi(f-h_d^2) \|_2\leq 2E_{{\alpha},d }+E_{{\alpha},d }^2. \]
	
	Let $e \coloneq f-|h_d|^2=f-h_d^2$, then clearly $e=e^*$, and
	\[ \| \Phi(e)\|_2^2 \geq   \| \Phi(e)^2\|_2=  \| \Phi(e^2)\|_2 =  \| \Phi(e^*e)\|_2.\]
	Because  $\| \Phi(e^*e)\|_2$ is the largest eigenvalue of $ \Phi(e^*e)$, and 
	\[C_{e^*e}(1_G)=\sum_{\rho \in \irr{G}} \frac{n_\rho}{{\card{G}}}\tr(\rho(e^*e))\]
	is the weighted mean value of \cyan{eigenvalues} of $ \Phi(e^*e)$, we have
	$\| \Phi(e)\|_2^2 \geq   \| \Phi(e^*e)\|  \geq C_{e^*e}(1_G)$.
	Since 
	\[e^*e=\sum_{x\in G}\sum_{y\in G}\overline{C_{e}(x)}C_{e}(y)x^{-1}y,\]
	we have 
	\[C_{e^*e}(1_G)=\sum_{x\in G}|C_{e}(x)|^2\leq E_{{\alpha},d }^4 + 4E_{{\alpha},d }^3 + 4E_{{\alpha},d }^2.\]
	The proof is complete by invoking Proposition~\ref{coro:appr-rate-sqrt}.
\end{proof}

Having established the error estimate for the thresholding hierarchy, \rred{the following result provides a sufficient hierarchy level that guarantees the existence of an approximately feasible primal-dual pair satisfying the stopping criteria commonly used by SDP solvers.} To this end, we first rewrite \eqref{Prob_feas_sparse_SDP_Gram} as 
\begin{eqnarray}\label{SDp:primal-of-vdQvd}
	\operatorname{minimize} \quad  &&\langle 0, Q \rangle \\
	\nonumber \operatorname{subject~to} \quad  &&Q \in \mathbb{C}^{S_d \times S_d},\; Q^\h = Q,\; Q \succeq 0, \\
	\nonumber    &&\mathcal{A}_d(Q)=C_f.
\end{eqnarray}
Here $\mathcal{A}_d$ is the linear operator corresponding to linear constraints in   \eqref{Prob_feas_sparse_SDP_Gram}. Correspondingly, the dual problem of \eqref{SDp:primal-of-vdQvd} is
\begin{eqnarray}\label{SDp:dual-of-vdQvd}
	\operatorname{maximize}  \quad  &&\langle C_f, y \rangle \\
	\nonumber \operatorname{subject~to} \quad && y\in  \mathbb{C}^{\cyan{S_d^{-1} \cdot S_d}}, \\
	\nonumber &&-\mathcal{A}_d^*(y) \succeq 0,
\end{eqnarray}
here $S_d^{-1}\cdot S_d=\{a^{-1}b: a,b \in S_d\}$, and $\mathcal{A}_d^*$ is the adjoint operator of $\mathcal{A}_d$. We recall from \cite{doi:10.1137/1038003} that when Problem~\eqref{SDp:primal-of-vdQvd} is strictly feasible, a pair $(Q,y)$ is primal-dual optimal if  it satisfies the KKT conditions:
\begin{equation}\label{eq:kkt:sat}
	\mathcal{A}_d(Q)=C_f,~Q\succeq 0,\quad 
	-\mathcal{A}_d^*(y) \succeq 0,\quad 
	\left< \mathcal{A}_d(Q),y\right>=\left< Q,\mathcal{A}_d^*(y)\right> =0.
\end{equation} 
In practice, numerical errors are unavoidable. Consequently, KKT conditions~\eqref{eq:kkt:sat} are discussed with respect to a numerical tolerance \cyan{\cite{lasserre2019sdp,doi:10.1137/15M1041924}}. Along this direction, we have the lemma that follows.  

\begin{lemma}\label{coro:numerical-feasible}
	Let $f \in \CG$. If the numerical tolerance $\eta $ satisfies that  \red{$\eta>2\sqrt{\varrho_d}$}, then there exists \cyan{a} matrix $\tilde Q\in\mathbb{C}^{S_d \times S_d}$ and \cyan{a} vector $\tilde y\in\mathbb{C}^{\cyan{S_d^{-1} \cdot S_d}}$ for the Problem~\eqref{Prob_feas_sparse_SDP_Gram} (or equivalently the Problem~\eqref{SDp:primal-of-vdQvd} and \eqref{SDp:dual-of-vdQvd}) such that
	\begin{eqnarray}
		\label{cond-num-fes1}	\|\mathcal{A}_d(\tilde Q)-C_f \|_2 \leq \eta, \label{eq:num:sat-1}	\\
		|\left<\mathcal{A}_d^*(\tilde y),\tilde Q\right> |\leq \eta,\\
		|\left<C_f,\tilde y\right> |\leq \eta,\\
		\label{cond-num-fes4}	\tilde Q \succeq 0,~-\mathcal{A}^*_d( \tilde y) \succeq 0.  \label{eq:num:sat-4}
	\end{eqnarray}
\end{lemma}
\begin{proof}
By the definition of $\varrho_d$, there exists $\tilde Q \in \mathbb{C}^{S_d \times S_d}$ such that
	\[\tilde Q \succeq 0,  \; \| \mathcal{A}_d(\tilde Q)-C_f\|_2 \leq 2\sqrt{\varrho_d} .\]
	Let  $\tilde y=0$, then $-\mathcal{A}_d^*(\tilde y) =0$, and therefore
	\[\left|\left<\mathcal{A}_d^*(\tilde y),\tilde Q\right> \right|= 	|\left<C_f,\tilde y\right> |=0.\]
	\rred{Clearly} $\tilde Q \succeq 0$, thus $(\tilde Q, \tilde y)$ satisfies \eqref{eq:num:sat-1}--\eqref{eq:num:sat-4}.
\end{proof} 

In fact, conditions \eqref{eq:num:sat-1}--\eqref{eq:num:sat-4} usually serve as stopping criteria in SDP solvers \cite{sun2020sdpnal+}. Moreover, this result can be directly applied to solvers, such as \cite{Garstka_2021}, that use the $\ell_\infty$-norm condition (i.e. $\|\mathcal{A}(Q) - b\|_\infty<\eta$) as a stopping criterion. Indeed, the inequality $\|\mathcal{A}(Q) - b\|_\infty \leq \| \mathcal{A}(Q) - b\|_2$ always holds. Therefore, any approximate solution satisfying the $\ell_2$-norm accuracy requirement automatically satisfies the corresponding $\ell_\infty$-norm accuracy requirement.

Next, we focus on $G = \red{\mathbb{Z}_2^n}$. In the fields of combinatorial optimization and theoretical computer science, many problems, including MAX-CUT and MAX-SAT, can be modeled as optimization problems over $\red{\mathbb{Z}_2^n}$.
\begin{corollary}\label{coro:binary:tolerance}
	Given $\alpha \in (0,1)$, there is a constant $c_\alpha$ such that for any $\eta > 0$ and any $f \in \mathbb{R}[x_1,\dots, x_n]_k$ with $\alpha \le f \le 1$ over the hypercube $\{-1,1\}^n$, 
	\cyan{the SDP feasibility problem corresponding to the $d$-th order Lasserre relaxation for certifying the nonnegativity of $f$ has a primal-dual pair satisfying}  \eqref{eq:num:sat-1}--\eqref{eq:num:sat-4}, as long as 
    \[\rred{ d\geq k\left\lceil \frac{   \log(9c_\alpha)-\log(\eta)}{\log(1+{\alpha^2})-\log(1-{\alpha^2})}\right\rceil.}\]
	\end{corollary}
\begin{proof}
	Let $c_\alpha$ be the constant in Theorem~\ref{thm:SDP_rate}. By assumption and the proof of Theorem~\ref{thm:SDP_rate}, we have $\sqrt{\varrho_d} < 9E_{\alpha,d}$. Then the conclusion follows immediately from Lemma~\ref{coro:numerical-feasible}.
\end{proof}

\subsection{Convergence rate of $\delta_d$}\label{subsec:convergence analysis2}

In this subsection, we provide an \rred{upper bound for} $\delta_d > 0$ such that $f+ \delta_d$  has an SOHS with support $S_d$. For any $f \in \CG$, we define the $\cyan{\ell_1}$-norm of $C_{f}$ as 
\[\|C_f\|_{1}=\sum_{x\in G} \left|C_{f}(x)\right|.\]
Moreover, for a subset $S\subseteq G$, we denote $S^{-1} \coloneqq \{x^{-1}:x \in S\}$.
\begin{lemma}\label{lem:ell_1-for-HSOS}
	Let $f \in \CG$ with $f=f^*$. For any set $S\subseteq G$ with 
	\[\operatorname{supp}(f) \subseteq S^{-1} S \coloneq \{x^{-1}y:~x,y \in S \}, \]
	$f+\|C_{f}\|_1$ has an SOHS with support $S$.
\end{lemma}
\begin{proof}
	We partition the set $\rred{S^{-1}  S}$ into the union of two disjoint sets $S_0$ and $\rred{(S^{-1}S)}\setminus  S_0$, where 
	\[S_0=\{x \in S^{-1} S: x=x^{-1}\},~\rred{(S^{-1}S)}\setminus  S_0=\{x \in S^{-1} \cdot S: x\neq x^{-1}\}.\]
	We further partition $\rred{(S^{-1}S)}\setminus  S_0$ into two disjoint sets $\rred{(S^{-1}S)}\setminus  S_0=S_1 \sqcup S_{-1}$ with  $S_{1}^{-1}=S_{-1}$. Since $f=f^*$, we have $\overline{C_{f}(x)}=C_{f}(x^{-1})$ for each $x \in G$. Thus,  
	\[f=\sum_{x \in S_0}C_{f}(x)x+\sum_{x \in S_1}C_{f}(x)x+\sum_{x \in S_{-1}}C_{f}(x)x=\sum_{x \in S_0}C_{f}(x)x+\sum_{x \in S_1}\left(C_{f}(x)x+\overline{C_{f}(x)}x^{-1}\right).\]
	For $x=a^{-1}b \in S_0$ with $a,b \in S$ and $C_f(x)\neq 0$, we have $C_{f}(x)=\overline{C_{f}(x)}$ and
	\[\left| C_{f}(x) \right|+ C_{f}(x)x=\frac{\left| C_{f}(x) \right|}{2} \left( a+\frac{  C_{f}(x)  }{\left| C_{f}(x) \right|} b \right)^*\left( a+\frac{  C_{f}(x)  }{\left| C_{f}(x) \right|} b \right).\]
	For $x=a^{-1}b \in S_1$ with $a,b \in S$ and $C_f(x)\neq 0$, 
	\[2\left| C_{f}(x) \right|+ C_{f}(x)x+\overline{C_{f}(x)}x^{-1}={\left| C_{f}(x) \right|} \left( a+\frac{  C_{f}(x)  }{\left| C_{f}(x) \right|} b \right)^*\left( a+\frac{  C_{f}(x)  }{\left| C_{f}(x) \right|} b \right). 
	\]
	Therefore, $f + \lVert C_f \rVert_1$ has an SOHS supported on $S$.
\end{proof}

\begin{lemma}\label{thm:sqrtf_restrict}
	Let $f \in \CG$ and $S \subseteq G$ be such that $f\succeq 0$, $S=S^{-1}$ and $\operatorname{supp}(f) \subseteq S$. We define the function $C_{\sqrt{f}}|_{S}$ on $G$ by 
	\begin{equation*}
		C_{\sqrt{f}}|_{S}(x)=
		\begin{cases}
			C_{\sqrt{f}}(x) \quad &x \in S, \\
			0 \quad &x\notin S.
		\end{cases} 
	\end{equation*}
	Then there is a constant $M \leq 2 \|C_{\sqrt{f}}\|_{1}  \|C_{\sqrt{f}}-C_{\sqrt{f}}|_{S}\|_{1}$ such that $f+M$ has an SOHS with support $S$.
\end{lemma}
\begin{proof}
	Since $f\succeq 0$, we must have $\sqrt{f}=\sqrt{f}^*$. Denote 
	\[
	\sqrt{f}|_{S}\coloneq \sum_{x \in G}C_{\sqrt{f}}|_{S}(x)x,\quad g\coloneq f-(\sqrt{f}|_{S})^2.
	\]
	Then clearly $\sqrt{f}|_{S}=(\sqrt{f}|_{S})^*$,		$\operatorname{supp}(g) \subseteq S^{-1}S$ and
	\[ M\coloneqq \|C_g\|_{1}\leq \left\|C_{\sqrt{f}}+C_{\sqrt{f}}|_{S} \right\|_{1} \left\| C_{\sqrt{f}}-C_{\sqrt{f}}|_{S} \right\|_{1}.\]
	By construction, $f + M = g+M +(\sqrt{f}|_{S})^2$ has an SOHS with support $S$.
\end{proof}

The following theorem is concerned with an upper bound of $\delta_d$. 
\begin{theorem}[Convergence rate of $\delta_d$]\label{corr:error-bound-general}
	Given $\alpha \in [0,1)$, there is a constant $c_\alpha$ such that for any  \rred{$f=f^* \in \CG$} with $\Idg \succeq f\succeq \alpha \Idg$, we have 
	\[
	\delta_{d} \le \begin{cases}
		\frac{c_\alpha}{\sqrt{2d-1}} \left( \frac{1-\alpha^2}{1 + \alpha^2} \right)^d \beta_{f,G} \quad &\text{if $0 < \alpha < 1$}, \\
		\frac{c_0}{d} \beta_{f,G} \quad &\text{if $\alpha =0$},
	\end{cases}
	\]
	where for every \rred{ $d \in \mathbb{N}$, $\beta_{f,G} \coloneqq 
	\cyan{\min} \left\lbrace
	2\|C_{\sqrt{f}}\|_{1} \left( \card{G} -\card{S_d}\right) , 3 { \card{\operatorname{supp}(f)}^{d}}
	\right\rbrace$ }.
\end{theorem}

\begin{proof}
	\rred{The case $f=0$ is trivial.}  Assume $\beta_{f,G} = 2\|C_{\sqrt{f}}\|_{1} \left( \card{G} -\card{S_d}\right)$. According to Proposition~\ref{prop:S_d-thresholdingd}, $x\notin S_d$ implies $|C_{\sqrt{f}}(x)|<E_{\alpha,d}$, from which we derive 
	\begin{eqnarray*}
		\left\| C_{\sqrt{f}} - C_{\sqrt{f}}|_{S_d} \right\|_{1} = \sum_{x \in G} \left| C_{\sqrt{f}}(x)-C_{\sqrt{f}}|_{S_d} (x)\right| 
		&=& \sum_{x \notin S_d} \left|C_{\sqrt{f}}(x)\right|  \\
		&\leq&  \left(\card{G}-\card{S_d}\right)\left( \max_{x \notin S_d} \left|C_{\sqrt{f}}(x)\right| \right) \\ 
		& \leq &\left(\card{G}-\card{S_d}\right) E_{\alpha,d}.
	\end{eqnarray*}
	By Lemma~\ref{thm:sqrtf_restrict}, $f+2\|C_{\sqrt{f}}\|_{1}\left(\card{G}-\card{S_d}\right) E_{\alpha,d}$ has an SOHS on $S_d$.

	Next, we suppose $\beta_{f,G} = 3 { \card{\operatorname{supp}(f)}^{d}}$ and denote $e_d \coloneqq f-\cyan{p_{\alpha,d}(f)^2}$. Theorem~\ref{thm:SDP_rate} implies
	\[\sum_{x\in G}|C_{e_d}(x)|^2\leq (2E_{\alpha,d}+E_{\alpha,d}^2)^2 \leq (3E_{\alpha,d})^2.\]
	From $f=f^*$, we derive $p_{\alpha,d}(f)^*=p_{\alpha,d}(f)$. Thus, $S_d=S_d^{-1}$ and
	\[\operatorname{supp}(e_d)\rred{\subseteq} \operatorname{supp}(f) \cup \operatorname{supp}(p_{\alpha,d}(f)^2) \subseteq  S^{-1}_{d}\cdot S_d .\]
	Since $f$ is nonzero and $f$ is nonnegative, we have $1_G \in \operatorname{supp}(f)$. We consider
	\begin{equation}\label{eq:Fd}
		F_d\coloneq \{\Pi_{i=1}^r x_i: x_1,\dots, x_r\in \operatorname{supp}(f),\; r \leq d \}.
	\end{equation}
	Then $S^{-1}_{d} \cdot S_d\subseteq F_{2d}$. By  the fact that $1_G \in \operatorname{supp}(f)$, we may conclude that \[\card{\operatorname{supp}(e_{d})}\leq \card{S_d^{-1} \cdot S_d} \leq   \card{F_{2d}}\leq \card{\operatorname{supp}(f)}^{2d}.\]
	This implies 
	\[\|C_{e_d}\|_{1}=\sum_{x\in \operatorname{supp}(e_{d})}|C_{e_d}(x)| \leq \sqrt{\card{\operatorname{supp}(e_{d})}} \sqrt{\sum_{x\in \operatorname{supp}(e_{d})}|C_{e_d}(x)|^2} \leq  3  \card{\operatorname{supp}(f)}^{d}E_{\alpha,d}.\]
	By Lemma~\ref{lem:ell_1-for-HSOS}, $\|C_{e_d}\|_{1}+e_d$ has an SOHS with support $S_d$.
\end{proof}

\begin{remark}
	For each $f\in \mathbb{C}[G]$ with $\Idg \succeq f \succeq 0$, we have $C_{f}(1_G)=\sum_{x \in G}|C_{\sqrt{f}}(x)|^2 \leq 1$. This implies $\|C_{\sqrt{f}}\|_{1}\leq \sqrt{\card{G}}$.
\end{remark}

It is possible to improve the upper bound of $\delta_d$ in Theorem~\ref{corr:error-bound-general} for special groups.  

\begin{corollary}\label{corr:error-bound-ZN}
	Let $G$ be the cyclic group $\mathbb{Z}_N$ generated by $y$ with $y^N = 1$, or the dihedral group $D_{2N}$ generated by $\sigma$ and $\tau$ with $\sigma^N=\tau^2=(\sigma\tau)^2= 1_G$.
	Given $\alpha \in [0,1)$, there is a constant $c_\alpha$ such that for each $f\in \mathbb{C}[G]$ with $\Idg \succeq f \succeq \alpha \Idg$, we have 
	\[
	\delta_{d} \le \begin{cases}
		\frac{c_\alpha}{\sqrt{2d-1}} \left( \frac{1-\alpha^2}{1 + \alpha^2} \right)^d \cyan{\omega_{d,f}}  \quad &\text{if $0 < \alpha < 1$}, \\
		\frac{c_0}{d}\cyan{\omega_{d,f}} \quad &\text{if $\alpha =0$}.
	\end{cases}
	\]
	Here for $G = \mathbb{Z}_N$, we define $\cyan{\omega_{d,f}} \coloneqq 3\sqrt{4dk+1}$, where $k$ is the smallest positive integer such that \[
	\operatorname{supp}(f) \subseteq \{y^r:0 \leq  r\leq k \} \cup \{y^{N-r}: 0 \leq r\leq k \}.\]
	For $G = D_{2N}$, we define $\cyan{\omega_{d,f}} \coloneqq 3\sqrt{8dk+1}$, where $k$ is the smallest positive integer such that \[
	\operatorname{supp}(f) \subseteq \{\sigma^r: 0 \leq  r\leq k \} \cup \{\sigma^{N-r}: 0 \leq  r\leq k \} \cup  \{\sigma^r\tau:0 \leq r\leq k \} \cup \{\sigma^{N-r}\tau: 0 \leq  r\leq k \}.
	\]

\end{corollary}
\begin{proof}
	Let $F_{2d}$ be the set defined by \eqref{eq:Fd}. According to the proof of Theorem~\ref{corr:error-bound-general}, it suffices to improve the estimate of $|F_{2d}|$. We notice that for $G = \mathbb{Z}_N$, it holds that 
\rred{\[{F_{2d}} \subseteq \{y^r: 0 \leq  r\leq 2dk \} \cup \{y^{N-r}: 0 \leq  r\leq 2dk  \} ,\]
}
which implies $\card{F_{2d}}\leq 4dk+1$.
	
	For $G = D_{2N}$, we observe that
	\[
	\sigma^a \sigma^b = \sigma^{a+b},\quad \sigma^a \tau \sigma^b = \sigma^{a-b}\tau,\quad \sigma^a  \sigma^b\tau = \sigma^{a+b}\tau,\quad \sigma^a \tau \sigma^b\tau = \sigma^{a-b}.
	\]
	for all $a,b \in \mathbb{N}$. This leads to $\card{F_{2d}}\leq 8dk+1$ as
	\[F_{2d} \subseteq \{\sigma^r: 0 \leq  r\leq 2dk \} \cup \{\sigma^{N-r}:0 \leq   r\leq 2dk \} \cup  \{\sigma^r\tau: 0 \leq  r\leq 2dk \} \cup \{\sigma^{N-r}\tau: 0 \leq  r\leq 2dk \}. \qedhere \]
\end{proof}

We notice that the upper bound of $\delta_d$ in Theorem~\ref{corr:error-bound-general} depends on $|G|$. In contrast, the bound in Corollary~\ref{corr:error-bound-ZN} is independent of $|G|$ when $G$ is cyclic or dihedral. For comparison, in \cite{MR3555385}, the authors investigated \emph{exact} FSOS for functions on $\mathbb{Z}_N$ \cyan{whose supports are} the same as \rred{those in} Corollary~\ref{corr:error-bound-ZN}. They obtained a nice result, showing that such $f$ always admits \emph{exact} FSOS with sparsity at most $3k\log(N/k)$.  Our results show that if we allow arbitrarily small \emph{perturbations} \cyan{$\varepsilon \rred{>} 0$}, then $f+\varepsilon$ will have an SOHS with sparsity \emph{independent} of $N$.

These theoretical results explain the phenomenon observed in earlier experiments (cf.~\cite[Section~4.1.1]{yang2022computing}). There, the approximate FSOS sparsity of positive functions on $\mathbb Z_N$ was found to be nearly constant under fixed degree bounds as $N$ increases. For completeness, we provide additional numerical examples below. For $N \in \{10^s: 3 \le s \le 7\}$ and $\alpha \in \{10^{-t}: 1\le t \le 4\}$, we generate $100$ random elements $f \in \mathbb{C}[\widehat{\mathbb{Z}_N}]$ supported on  
\[
\{\chi_{l}(x) = \exp(2\pi i l x / N) : -5 \leq l \leq 5\},
\]
subject to $\Idg \succeq f \succeq \alpha \Idg$. 
Therefore the sparsity of each  randomly generated element is at most $11$. Moreover, for any $k \in\mathbb{N}$, the level-$k$ support in our hierarchy satisfies
\[  S_k\subseteq \{\chi_l:-5k\le l\le 5k\}, \] 
whose cardinality is at most $10k+1$. In our numerical experiments, the hierarchy usually terminates at level $k \leq 3 $, and the largest termination level observed is $6$. Thus the supports of SOHS obtained in these examples are sparse compared to the \rred{group} order $N$. For each $f$, we increase $d$, starting from $1$, until $f$ admits an SOHS decomposition supported on $S_d$. Then we record the sparsity of the resulting SOHS representation. To verify whether $f$ has an SOHS decomposition over $S_d$, we solve a feasibility SDP using the COSMO solver~\cite{Garstka_2021}. The absolute and relative tolerances of the solver are set to $10^{-6}$ (compared to the default value of $10^{-5}$ in COSMO), thus ensuring greater precision. 
   For varying values of $N$ and $\alpha$, the mean \rred{sparsity over 100 instances} is recorded in Table~\ref{tab:Numerical}. The fourth column of Table~\ref{tab:Numerical} is the theoretical upper bound on FSOS sparsity established in~\cite{MR3555385}. For better illustration, we also plot the mean value of the sparsity versus $N$ in Figure~\ref{fig:line_chart_log_x}, for $\alpha = 10^{-1},10^{-2},10^{-3}$ and $10^{-4}$. The code used to produce Table 3 and Figure 1 \rred{is} available at \url{https://github.com/jty-AMSS/SparseGroupSOHS}. 
   \begin{table}[ht]
	\centering
	\caption{Numerical experiments for $100$ random $f \in \mathbb{C}[\widehat{\mathbb{Z}_N}]$ with $\supp(f) \subseteq \{\chi_l:-5\le l\le 5\}$ and $\alpha\leq f\leq 1$.}
	\label{tab:Numerical}
	\begin{tabular}{|c|c|c|c|}
		\hline
		$N$ & $\alpha$ & {average approximate SOHS sparsity} & bounds \cite{MR3555385} \\
		\hline
		$10^3$ & $10^{-4}$ & 14.2 & 120 \\ \hline
		$10^3$ & $10^{-3}$ & 12.6 & 120 \\ \hline
		$10^3$ & $10^{-2}$ & 11.0 & 120 \\ \hline
		$10^3$ & $10^{-1}$ & 11.0 & 120 \\ \hline
		
		$10^4$ & $10^{-4}$ & 14.2 & 165 \\ \hline
		$10^4$ & $10^{-3}$ & 12.2 & 165 \\ \hline
		$10^4$ & $10^{-2}$ & 11.0 & 165 \\ \hline
		$10^4$ & $10^{-1}$ & 11.0 & 165 \\ \hline
		
		$10^5$ & $10^{-4}$ & 15.3 & 225 \\ \hline
		$10^5$ & $10^{-3}$ & 12.1 & 225 \\ \hline
		$10^5$ & $10^{-2}$ & 11.0 & 225 \\ \hline
		$10^5$ & $10^{-1}$ & 11.0 & 225 \\ \hline
		
		$10^6$ & $10^{-4}$ & 15.2 & 270 \\ \hline
		$10^6$ & $10^{-3}$ & 11.9 & 270 \\ \hline
		$10^6$ & $10^{-2}$ & 11.0 & 270 \\ \hline
		$10^6$ & $10^{-1}$ & 11.0 & 270 \\ \hline
		
		$10^7$ & $10^{-4}$ & 14.8 & 315 \\ \hline
		$10^7$ & $10^{-3}$ & 13.1 & 315 \\ \hline
		$10^7$ & $10^{-2}$ & 11.0 & 315 \\ \hline
		$10^7$ & $10^{-1}$ & 11.0 & 315 \\ \hline
	\end{tabular}
\end{table}
 
\pgfplotsset{
    compat=1.17,  %
    every axis plot/.append style={
        line width=1.2pt,  %
        mark size=4pt,     %
        mark options={fill=white}  %
    },
    every axis/.append style={
        axis lines=left,    %
        tick label style={font=\small},  %
        label style={font=\small},       %
        legend style={font=\small},      %
        grid=major,         %
        grid style={dashed, gray!30}     %
    }
}

\begin{figure}[htbp] %
    \centering  %
    \begin{tikzpicture}
        \begin{axis}[
            xlabel={N},  
            ylabel={average approximate SOHS sparsity},
            xmode=log,  %
            xmin=1000, xmax=1e7,  %
           xticklabels={$10^3$, $10^4$, $10^5$, $10^6$, $10^7$},  %
            xtick={1e3, 1e4, 1e5, 1e6, 1e7},  %
            ymin=10.0, ymax=16.0,  %
            ytick={10.0, 11.0, 12.0, 13.0, 14.0, 15.0,16.0},  %
            legend style={
                at={(1.05, 0.5)},   
                anchor=north west, 
                draw=none         
            },
            legend entries={$\alpha=10^{-4}$,$\alpha=10^{-3}$, $\alpha=10^{-2}$ ,$\alpha=10^{-1}$,}  
        ]

\addplot[purple ] coordinates {
    (1e3, 14.2)
    (1e4, 14.2)
    (1e5, 15.3)
    (1e6, 15.2)
    (1e7, 14.8)
};

\addplot[blue] coordinates {
    (1e3, 12.6)
    (1e4, 12.2)
    (1e5, 12.1)
    (1e6, 11.9)
    (1e7, 13.1)
};

\addplot[red] coordinates {
    (1e3, 11.0)
    (1e4, 11.0)
    (1e5, 11.0)
    (1e6, 11.0)
    (1e7, 11.0)
};

\addplot[green!50!black,dashed] coordinates {
    (1e3, 11.0)
    (1e4, 11.0)
    (1e5, 11.0)
    (1e6, 11.0)
    (1e7, 11.0)
};

        \end{axis}
    \end{tikzpicture}
    \caption{
Numerical experiments for $100$ random $f \in \mathbb{C}[\widehat{\mathbb{Z}_N}]$ supported on $\{\chi_l:-5\le l\le 5\}$ satisfying $\alpha\leq f\leq 1$, with $N \in \{10^s: 3 \le s \le 7\}$ and $\alpha \in \{10^{-t}: 1\le t \le 4\}$.
  }
    
    \label{fig:line_chart_log_x}  
\end{figure}

It should be emphasized that \rred{the numerical results are computed by numerical solvers and therefore unavoidably contain numerical errors,} whereas the analysis in \cite{MR3555385} provided a theoretical upper bound for \emph{exact} FSOS.
Therefore, we cannot attribute the gap  between experimental results and the theoretical upper bound to the insufficient tightness of the theoretical upper bound in \cite{MR3555385}.

\section*{Further discussion}
It is possible to extend the framework \cite{MR3555385,MR3629437,MR5048294} for sum-of-squares in finite abelian group algebras to the setting of finite non-abelian group algebras. To be more precise, we denote by $\operatorname{Cay}(G, \operatorname{supp}(f))$ the Cayley graph of $G$ with respect to $\operatorname{supp}(f)$. The same arguments as in \cite{MR3555385} show that one possible approach to constructing an exact sparse SOHS certificate of $f$ is to find a small chordal cover of $\operatorname{Cay}(G, \operatorname{supp}(f))$ and decompose it into maximal cliques. This approach is quite effective for finite abelian groups such as $\mathbb{Z}_N$ and $\mathbb{Z}_2^n$, for which small chordal covers can be explicitly constructed \cite{MR3555385,MR5048294}. For general graphs, however, finding a chordal cover with small maximal cliques is equivalent to the treewidth problem, which is NP-hard \cite{ACP87}. Therefore, for a general finite non-abelian group $G$, we do not expect the exact sparse SOHS problem to be efficiently solvable by the Cayley graph method. Moreover, since this paper studies SOHS certificates with error rather than exact SOHS certificates, the aforementioned framework would require additional modifications. As far as we are aware, such an adaptation has not yet been developed even for finite abelian groups. 

\subsection*{Acknowledgment}
We are grateful to the editors and reviewers for their helpful comments, which have substantially improved the paper.

\bibliographystyle{plain}
\bibliography{optimization}
\end{document}